\documentclass[a4paper, 11pt]{article}
\usepackage{amsfonts,amsmath,amssymb,amscd,color,stmaryrd,latexsym}
\usepackage[all]{xy} 
\usepackage[T1]{fontenc}
\usepackage{hyperref}

\def\End{\mathop{\rm End}\nolimits}

\def\vol{\mathop{vol}\nolimits}
\def\ker{\mathop{\rm ker}\nolimits}
\def\im{\mathop{\rm im}\nolimits}

\def\diff{\mbox{\sl Diff}}

\def\note#1{\marginpar{\raggedright\if@twoside\ifodd\c@page\raggedleft\fi\fi\sf\scriptsize RMK: #1}}

\newcommand{\SU}{\mr{SU}}
\newcommand{\SO}{\mr{SO}}

\newcommand{\GL}{\mr{GL}}

\newcommand{\Gt}{\mr{G_2}}
\newcommand{\spin}{\mr{Spin}}

\newcommand{\id}{\text{Id}}

\newcommand{\mf}{\mathfrak}
\newcommand{\mr}{\mathrm}
\newcommand{\mb}{\mathbf}
\newcommand{\mc}{\mathcal}

\newcommand{\R}{\mathbb{R}}

\newcommand{\C}{\mathbb{C}}

\renewcommand{\iff}{if and only if }
\newcommand{\wrt}{with respect to }
\newcommand{\st}{such that }

\newcommand{\Id}{\mr{Id}}

\newcommand{\loc}{\stackrel{\mr{loc}}{=}}

\newcommand{\be}{\begin{equation}}
\newcommand{\ben}{\begin{equation}\nonumber}
\newcommand{\ee}{\end{equation}}

\newcommand{\grad}{\mr{grad}\,}
\newcommand{\Db}{\mbox{\DH}}

\newenvironment{rem}{\begin{trivlist}\item[]{\bf Remark:}\setlength{\parindent}{0pt}}{\end{trivlist}}
\newenvironment{ex}{\begin{trivlist}\item[]{\bf Example:}\setlength{\parindent}{0pt}}{\end{trivlist}}
\newenvironment{prf}{\begin{trivlist}\item[]{\bf Proof:}}{\hfill$\blacksquare$\end{trivlist}}

\newenvironment{thm*}{\begin{trivlist}\item[]{\bf Theorem.}}{\end{trivlist}}
\newenvironment{dfn*}{\begin{trivlist}\item[]{\bf Definition.}}{\end{trivlist}}
\newenvironment{prp*}{\begin{trivlist}\item[]{\bf Proposition.}}{\end{trivlist}}
\newenvironment{lem*}{\begin{trivlist}\item[]{\bf Lemma.}}{\end{trivlist}}
\newenvironment{cor*}{\begin{trivlist}\item[]{\bf Corollary.}}{\end{trivlist}}

\newtheorem{thm}{Theorem}[section]
\newtheorem{dfn}[thm]{Definition}
\newtheorem{prp}[thm]{Proposition}
\newtheorem{lem}[thm]{Lemma}
\newtheorem{cor}[thm]{Corollary}

\begin{document}
\title{A heat flow for special metrics}
\author{Hartmut Wei{\ss} and Frederik Witt}
\date{\today}
\maketitle

\centerline{\textbf{Abstract}}
On the space of positive $3$--forms on a seven--manifold, we study a natural functional whose critical points induce metrics with holonomy contained in $\Gt$. We prove short--time existence and uniqueness for its negative gradient flow. Furthermore, we show that the flow exists for all times and converges modulo diffeomorphisms to some critical point for any initial condition sufficiently $C^\infty$--close to a critical point.

\bigskip

\textsc{MSC 2000:} 53C44 (35G25, 35K55, 53C25)

\smallskip

\textsc{Keywords:} $\Gt$--manifolds,  geometric evolution equations, quasilinear parabolic equations
%
%
%
%
%
\section{Introduction}
A central problem in Riemannian geometry is the construction of metrics with prescribed properties of the Ricci tensor. In this context the metrics we refer to as {\em special} are particularly interesting. These metrics are induced by a differential form $\Omega$ of special algebraic type subject to the non--linear harmonic equation
\be\label{basiceq}
d\Omega=0,\quad d\star_\Omega\!\Omega=0,
\ee
where $\star_\Omega$ is the Hodge operator associated with the induced metric $g_\Omega$. A prototypical example are so--called $\Gt$--metrics. A theorem of Fern\'andez and Gray~\cite{fegr82} asserts~\eqref{basiceq} to hold \iff the holonomy of $g_\Omega$ is contained in $\Gt$. This, in turn, implies $g_\Omega$ to be Ricci flat~\cite{bo66}. Other metrics of this type are $\spin(7)$--metrics (which are also Ricci flat) or quaternionic K\"ahler metrics (which are Einstein), but more exotic examples such as $\mathrm{PSU}(3)$--metrics~\cite{hi01} (satisfying a less standard condition on the Ricci tensor)  also fit into this setting.

\bigskip

In the $\Gt$--case we consider {\em positive} $3$--forms $\Omega$ over an oriented $7$--manifold which are special in as far they give rise to a complementary $4$--form $\Theta(\Omega)$ so that the volume form $\vol_\Omega:=\Omega\wedge\Theta(\Omega)/7$ induces the chosen orientation. In fact, once the metric has been constructed from $\Omega$ one has $\Theta(\Omega)=\star_\Omega\Omega$. If $M$ is compact,~\eqref{basiceq} is equivalent to the non--linear Laplace equation $\Delta_\Omega\Omega=0$, where $\Delta_\Omega$ is the Hodge Laplacian associated with $g_\Omega$. In this sense, $\Omega$ is ``self--harmonic'', but we shall stick to the usual $\Gt$--jargon and refer to positive forms satisfying~\eqref{basiceq} as {\em torsion--free}. For closed $\Omega$, Hitchin~\cite{hi01} interpreted the second condition as the Euler--Lagrange equation for the functional on positive $3$--forms $\Omega\mapsto\int_M\vol_\Omega$ restricted to the cohomology class $[\Omega]$. Existence of critical points, however, is a delicate issue. Since Joyce's seminal work~\cite{jo00} we know non--trivial compact holonomy $\Gt$--manifolds to exist, but a Yau--Aubin type theorem which guarantees a priori existence is yet missing.

\bigskip

A natural idea for proving existence of $\Gt$--metrics is to look for a geometric evolution equation on the space of positive $3$--forms, which could evolve forms towards a torsion--free $\Omega$, cf.~\cite{ka09} for a rather general perspective. A first candidate for a flow equation has been proposed by Bryant in~\cite{br06}, namely the ``Laplacian flow''
\ben
\frac{\partial}{\partial t}\Omega=\Delta_\Omega\Omega.
\ee
Restricted to closed positive $3$--forms, we can think of this flow as the ($L^2$--)gradient flow of Hitchin's functional. However, as we are going to show, the resulting flow equation is not even weakly parabolic so that standard parabolic theory does not apply directly (though Bryant  and Xu, resp.~Xu and Ye have managed to establish short--time existence for the Laplacian flow for {\em closed} initial conditions, cf.~\cite{brxu11,xuye09}). This is reminiscent of the Einstein--Hilbert functional whose negative gradient is difficult to deal with on the same grounds, a fact which subsequently led to the definition of Ricci flow. We therefore consider the negative gradient flow of the ``Dirichlet energy'' functional
\ben
\mc{D}:\Omega\mapsto\frac{1}{2}\int_M\big(|d\Omega|_{g_\Omega}^2+|d\Theta(\Omega)|_{g_\Omega}^2\big)\vol_\Omega
\ee
whose critical points, as we will show, are precisely given by torsion--free forms. In principle, the definition of $\mc{D}$ makes sense for any special metric. The reason to focus on $\Gt$ is twofold. The set of positive $3$--forms $\Omega$ is an open subset of $\Omega^3(M)$, and $\Gt$ acts transitively on the sphere (as do all reduced holonomy groups of manifolds which are not locally symmetric). Both these features greatly simplify technicalities.

\bigskip

Our first result is short--time existence and uniqueness.

\begin{thm*}
Let $M$ be a closed, oriented $7$--manifold. Given a positive $3$--form $\Omega_0$ on $M$, there exists $\epsilon>0$ and a smooth family of positive $3$--forms $\Omega(t)$ for $t\in[0,\epsilon]$ \st
\ben
\frac{\partial}{\partial t}\Omega=-\grad\mc{D}(\Omega),\quad\Omega(0)=\Omega_0.
\ee  
Furthermore, for any two solutions $\Omega(t)$ and $\Omega'(t)$ we have $\Omega(t)=\Omega'(t)$ whenever defined. 
\end{thm*}

Hence we can speak of {\em the} Dirichet energy flow for some initial condition $\Omega_0$ defined on a maximal time--interval $[0,T_{max})$, $0<T_{max}\leq\infty$. We emphasise that for the Dirichlet energy flow the initial condition does {\em not} have to be assumed to be closed. The proof is based on DeTurck's trick as introduced in~\cite{dt83}, namely to consider a geometric perturbation $\widetilde Q$ of the negative gradient of $\mc{D}$. Then $\widetilde Q$ is strongly elliptic and the standard theory of quasilinear parabolic equations applies. Of course, one cannot expect longtime existence and convergence in general as the existence of $\Gt$--metrics is topologically obstructed. For instance, if $M$ is compact then $\Omega$ cannot be exact, hence the third Betti number $b_3$ must be greater or equal to $1$. However, a meaningful negative gradient flow is certainly expected to exist for all times and to converge near a critical point. Indeed we will prove the subsequent stability result (cf.~Theorem~\ref{stability} and Corollary~\ref{convcor} for the precise statement).

\begin{thm*}
Let $M$ be a closed, oriented $7$--manifold and let $\bar\Omega$ be a torsion--free $\Gt$--form on $M$. For initial conditions sufficiently $C^\infty$--close to $\bar\Omega$ the Dirichlet energy flow exists for all times and converges modulo diffeomorphisms to a torsion--free $\Gt$--form.
\end{thm*}

The theorem resembles \v Se\v sum's corresponding stability result for Ricci flat metrics under Ricci flow~\cite{se04}. Nevertheless, two major differences occur: For stability of the Ricci flow one needs to assume that (a) the Ricci flat metric $\bar g$ (which corresponds to the torsion--free $\bar\Omega$ in our setting) is a smooth point of the moduli space of Ricci flat metrics (b) the linearisation at $\bar g$ of the quasilinear elliptic operator involved (corresponding to our $\widetilde Q$) is non--positive. Both conditions are difficult to check in practice unless one makes further assumptions such as special holonomy. In our situation, however, these assumptions are automatically satisfied. In particular, we can show that near $\bar\Omega$, the set $\widetilde Q^{-1}(0)$ provides a slice for the action of diffeomorphisms isotopic to the identity on the space of torsion--free $\Gt$--metrics. Consequently, the moduli space of torsion--free $\Gt$--forms is smooth. This has been proven previously by Joyce (see~\cite{jo00}), but our approach is rather different and based on the specific geometry of the zero sets of $\grad\mc{D}$ and $\widetilde Q$. Secondly, we note that on a technical level, the linearisation $L_{\bar\Omega}$ of $\widetilde{Q}_{\bar\Omega}$ at $\bar\Omega$ does not admit a Weitzenb\"ock formula as does the Lichnerowicz Laplacian which appears in the Ricci flow case. However, we can prove the G\r arding type inequality
\ben
\langle-L_{\bar\Omega}\dot\Omega,\dot\Omega\rangle_{L^2}\geq C\|\dot\Omega\|^2_{W^{1,2}}-\|\dot\Omega\|^2_{L^2},
\ee
where $\dot\Omega\in\Omega^3(M)$ and $C$ is some positive constant. With this coercivity condition we can invoke a result of Lax--Milgram type in order to establish longtime existence following ideas developed in~\cite{hp99}. Finally, convergence follows from a careful analysis of the remainder term $R_{\bar\Omega}=\widetilde Q-L_{\bar\Omega}$. The flow is thus, in principle, capable of finding torsion--free $\Gt$--forms on a $7$--manifold if they exist.
%
%
%
%
%
\section{$\mb{G_2}$--structures}\label{g2s}
Here and in the sequel, we let $M$ be a connected, closed, oriented manifold of dimension $7$. We first recall some basic features of $\Gt$--geometry to fix notations. Good references are~\cite{br87} and Chapter 10 in~\cite{jo00}.

\bigskip

The group $\GL(7)$ acts on $\Lambda^3\R^{7*}$ and has an open orbit $\mc{O}$ diffeomorphic to $\GL(7)/\Gt$. In fact, $\mc{O}$ is an open cone, for $c\cdot\Id\in\GL(7)$, $c \neq0$ acts on any $3$--form by multiplication with $c^{-3}$. Since $\Gt$ is a subgroup of $\SO(7)$, any $\Omega\in\mc{O}$ induces an orientation and a Euclidean metric $g_{\Omega}$ on $\R^7$. We denote by $\Lambda^3_+$ those $\Omega\in\mc{O}$ which induce the standard orientation on $\R^7$ and refer to its elements as {\em positive} forms. They are acted on transitively by $\GL(7)_+$, the orientation--preserving linear isomorphisms of $\R^7$.

\bigskip

Let $\Omega_+^3(M)$ denote the open set of sections of $\Lambda^3_+M$, the fibre bundle associated with $\Lambda^3_+$. Then a section $\Omega\in\Omega_+^3(M)$ (which exists \iff the second Stiefel--Whitney class of $M$ vanishes) induces a reduction of the frame bundle to a principal $\Gt$--bundle. We also refer to the pair $(M,\Omega)$ as a {\em $\Gt$--structure}. Such a structure singles out a principal $\SO(7)$--bundle whose associated metric, Hodge star operator and Levi--Civita connection we denote by $g_\Omega$, $\star_\Omega$ and $\nabla^\Omega$ respectively. Locally, there exist so--called {\em $\Gt$--frames}, i.e.\ local frames $(e_1,\ldots,e_7)$ of $TM$ for which $\Omega$ has ``normal form''. In our convention, we can write $\Omega$ \wrt a $\Gt$--frame as
\ben
\Omega=e^{127}+e^{347}+e^{567}+e^{135}-e^{146}-e^{236}-e^{245}
\ee
with $e^{ijk}$ shorthand for $e^i\wedge e^j\wedge e^k$. Note that a $\Gt$--frame is orthonormal for the induced $\Gt$--metric $g_\Omega$.

\bigskip

The holonomy of $g_\Omega$ is contained in $\Gt$ (implying that $g_\Omega$ is Ricci--flat~\cite{bo66}) \iff the underlying $\Gt$--form $\Omega$ is parallel, i.e.\ $\nabla^{g_\Omega}\Omega=0$. In this case we shall say that the $\Gt$--structure is {\em torsion--free} while we call $(M,\Omega)$ a {\em holonomy $\Gt$--manifold} if the holonomy of $g_\Omega$ is actually equal to $\Gt$. (This convention is by no means universal in the literature.) A torsion--free $\Gt$--structure has holonomy $\Gt$ if the fundamental group $\pi_1(M)$ is finite ($M$ being compact). By a theorem of Fern\'andez and Gray~\cite{fegr82}, torsion--freeness is equivalent to $d\Omega=0$ and $\delta_\Omega\!\Omega=0$, where $\delta_\Omega=(-1)^p\star_\Omega \!d\,\star_\Omega$ is the induced codifferential on $p$--forms. We shall therefore refer to any such $\Omega$ as a {\em torsion--free} $\Gt$--form. The latter equation can be viewed as the Euler--Lagrange equation of a non--linear variational problem set up by Hitchin~\cite{hi01}. Consider the smooth $\GL(7)_+$--equivariant map
\ben
\phi:\Lambda^3_+\to\Lambda^7,\quad\Omega\mapsto\vol_\Omega:=\star_\Omega1=\tfrac{1}{7}\Omega\wedge\star_\Omega\Omega,
\ee
whose first derivative at $\Omega$ evaluated on $\dot\Omega\in\Lambda^3$ is
\be\label{phider}
D_\Omega\phi(\dot\Omega)=\frac{1}{3}\star_\Omega\!\Omega\wedge\dot\Omega.
\ee
Integrating $\phi$ gives the functional
\be\label{functional}
\mc{H}:\Omega_+^3(M)\to\R,\quad\Omega\mapsto\int_M\phi(\Omega).
\ee
In analogy with Hodge theory we can restrict $\mc{H}$ to a fixed cohomology class and ask for critical points. From~\eqref{phider} it follows that a closed $\Omega$ is a critical point in its cohomology class \iff $\delta_\Omega\!\Omega=0$~\cite{hi01}. In particular, $\Omega$ is torsion--free and thus harmonic \wrt its induced Laplacian $\Delta_\Omega=d\delta_\Omega+\delta_\Omega d$. In passing we note that modulo a constant, $\mc{H}(\Omega)$ can be viewed as the norm square of the Euler vector field $\Omega \mapsto \Omega$ on the ``prehilbert'' manifold $\Omega_+^3(M)$ with induced $L^2$--metric 
\ben
\langle\dot\Omega_1,\dot\Omega_2\rangle_{L^2_\Omega}:=\int_Mg_\Omega(\dot\Omega_1,\dot\Omega_2)\vol_\Omega=\int_M\dot\Omega_1\wedge\star_\Omega\dot\Omega_2,
\ee
for elements $\dot\Omega_1$, $\dot\Omega_2$ in the tangent space $T_\Omega\Omega_+^3(M)\cong\Omega^3(M)$. We will drop the reference to $\Omega$ whenever this can be safely done and simply write $\langle\cdot\,,\cdot\rangle_{L^2}$ and $g$. The associated norms are then denoted by $\|\cdot\|$ and $|\cdot|$ respectively. 
%
%
%
%
%
\section{Representation theory}\label{reptheory}
Next we recall some elements of $\Gt$--representation theory. Most of the material is standard (mainly taken from~\cite{br06}, \cite{chsa02} and \cite{ka10}) or follows from straightforward computations.

\bigskip

The group $\Gt$ acts irreducibly in its vector representation $\Lambda^1\cong\R^7$ (in presence of a metric, we tacitly identify vectors with their duals). This action extends to the exterior algebra in the standard fashion, though $\Lambda^p$, the $\Gt$--representation over $p$--forms, is no longer irreducible for $2\leq p\leq5$. More precisely, we have orthogonal decompositions
\ben
\Lambda^2=\Lambda^2_7\oplus\Lambda^2_{14},\quad\Lambda^3=\Lambda^3_1\oplus\Lambda^3_7\oplus\Lambda^3_{27},
\ee
where the subscript indicates the dimension of the module. We denote the corresponding components by $[\alpha^p]_q$. For $p=2,3$ they satisfy
\be\label{star2forms}
\begin{array}{ll}
\Lambda^2_7=\{\alpha\in\Lambda^2\,|\,\star_\Omega(\alpha\wedge\Omega)=2\alpha\},&\Lambda^2_{14}=\{\alpha\in\Lambda^2\,|\,\star_\Omega(\alpha\wedge\Omega)=-\alpha\},\\[5pt]
\Lambda^3_7=\{\star_\Omega(X\wedge\Omega)\,|\,X\in\Lambda^1\},&\Lambda^3_{27}=\{\alpha\in\Lambda^3\,|\,(\star_\Omega\Omega)\wedge\alpha=0,\,\Omega\wedge\alpha=0\}.
\end{array}
\ee
The Lie algebra of $\mf{g}_2$ sitting inside $\mf{so}(7)\cong\Lambda^2$ corresponds to $\Lambda^2_{14}$, while $\Lambda^3_1$ simply consists of multiples of $\Omega$. Note that by equivariance, $\star_\Omega$ induces isomorphisms $\Lambda^p_q\cong\Lambda^{7-p}_q$ from which an analogous decomposition of $\Lambda^4$ and $\Lambda^5$ follows. This and the characterisations~\eqref{star2forms} are obtained from a routine application of Schur's lemma. For illustration, we derive for $\eta\in\Lambda^2$ the identity
\be\label{7contrOm}
(\eta\llcorner\Omega)\llcorner\Omega=3[\eta]_7.
\ee
Here, $\llcorner$ denotes the extension of the metric contraction to $\llcorner:\Lambda^k V^*\otimes\Lambda^l V^*\to\Lambda^{l-k}V^*$, e.g.\ $e^{12}\llcorner e^{12345}=e^{345}$ etc. Now $\eta\mapsto\eta\llcorner\Omega$ is a $\Gt$--equivariant map taking values in the irreducible module $\Lambda^1=\Lambda^1_7$ so that by Schur $\Lambda^2_{14}\subset\ker\,\llcorner\Omega$, whence $\eta\llcorner\Omega=[\eta]_7\llcorner\Omega$. Therefore, the identity~\eqref{7contrOm} needs only to be checked for one nontrivial element in $\Lambda^2_7$ (again by Schur). Fixing a $\Gt$--frame as in the previous section, we find $e_1\llcorner\Omega=e_{27}+e_{35}-e_{46}\in\Lambda^2_7$, hence $(e_1\llcorner\Omega)\llcorner\Omega=3e_1$. In the same vein, we can prove
\begin{eqnarray}\label{morecontr}
\big(\star_\Omega(\alpha\wedge\Omega)\big)\wedge\Omega & = & -4\star_\Omega\alpha,\nonumber\\
\big(\star_\Omega(\alpha\wedge\star_\Omega\Omega)\big)\wedge\star_\Omega\Omega & =& 3\star_\Omega\!\alpha,\\
\big(\star_\Omega(\alpha\wedge\star_\Omega\Omega)\big)\wedge\Omega & = & 2\alpha\wedge\star_\Omega\Omega.\nonumber  
\end{eqnarray}

\medskip

If the manifold $M$ is endowed with a $\Gt$--structure $\Omega$, all these decompositions and identities acquire global meaning. In particular we can speak of $\Omega^p_q$--forms, where $\Omega^p_q(M)=C^{\infty}(\Lambda^p_qT^*\!M)$ are smooth sections of the bundles with fibre $\Lambda^p_q$. As in the case for K\"ahler manifolds, this decomposition gives rise to $\Gt$--analogues of the Cauchy--Riemann operator, provided the $\Gt$--structure is torsion--free. The subsequent formul{\ae} were derived by Bryant and Harvey and can be found in~\cite{br06}. We briefly describe their results. First, we fix reference modules for the irreducible $\Gt$--representations occuring in the exterior algebra $\Lambda^*$, namely $\Omega_1=\Omega^0_1(M)$, $\Omega_7=\Omega^1_7(M)$, $\Omega_{14}=\Omega^2_{14}(M)$ and $\Omega_{27}=\Omega^3_{27}(M)$. Any form $\alpha\in\Omega^p(M)$ can be written in terms of $\Gt$--equivariant maps applied to elements of these reference modules. For instance $\dot\Omega=[\dot\Omega]_1\oplus[\dot\Omega]_7\oplus[\dot\Omega]_{27}\in\Omega^3(M)$ can be written as $\dot\Omega=f\Omega\oplus\star_\Omega(\alpha\wedge\Omega)\oplus\gamma$ for $f\in\Omega_1$, $\alpha\in\Omega_7$ and $\gamma\in\Omega_{27}$. There exist first order differential operators $d^p_q:\Omega_p\to\Omega_q$ such that the identities of Table~\ref{extderfor} hold.

\medskip

\begin{table}[htb]
  \quad\begin{tabular}{ l c l l l l}
      	$df$ & $=$ & & $\phantom{+}d^1_7f$ & &\\[2pt]
	$d(f\Omega)$ & $=$ & & $\phantom{+}d^1_7f\wedge\Omega$ & &\\[2pt]
	$d(f\!\star_\Omega\!\Omega)$ & $=$ & & $\phantom{+}d^1_7f\wedge\star\Omega$ & &\\[7pt]
	$d\alpha$ & $=$ & & $\phantom{+}\tfrac{1}{3}\star_\Omega(d^7_7\alpha\wedge\star_\Omega\Omega)$ & $+d^7_{14}\alpha$&\\[2pt]
	$d\star_\Omega(\alpha\wedge\star_\Omega\Omega)$ & $=$ & $-\tfrac{3}{7}d^7_1\alpha\cdot\Omega$ & $-\tfrac{1}{2}\star_\Omega(d^7_7\alpha\wedge\Omega)$ & & $+d^7_{27}\alpha$\\[2pt]
	$d\star_\Omega(\alpha\wedge\Omega)$ & $=$ & $\phantom{-}\tfrac{4}{7}d^7_1\alpha\cdot\star_\Omega\Omega$ &$+\frac{1}{2}d^7_7\alpha\wedge\Omega$ & &$+\star_\Omega d^7_{27}\alpha$\\[2pt]
	$d(\alpha\wedge\Omega)$ & $=$ & & $\tfrac{2}{3}d^7_7\alpha\wedge\star_\Omega\Omega$ & $-\star_\Omega d^7_{14}\alpha$ &\\[2pt]
	$d(\alpha\wedge\star_\Omega\Omega)$ & = & & $\phantom{+}\star_\Omega d^7_7\alpha$ & &\\[2pt]
	$d(\star_\Omega\alpha)$ & $=$ & $-d^7_1\alpha\cdot\vol_\Omega$\\[7pt]
	$d\beta$ & $=$ & & $\phantom{+}\tfrac{1}{4}\star_\Omega(d^{14}_7\beta\wedge\Omega)$ & &$+d^{14}_{27}\beta$\\[2pt]
	$d(\star_\Omega\beta)$ & $=$ & & $\phantom{+}\star_\Omega d^{14}_7\beta$ & &\\[7pt]
	$d\gamma$ & $=$ & & $\phantom{-}\tfrac{1}{4}d^{27}_7\gamma\wedge\Omega$ & &$+\star_\Omega d^{27}_{27}\gamma$\\[2pt]
	$d(\star_\Omega\gamma)$ & $=$ & & $-\tfrac{1}{3}d^{27}_7\gamma\wedge\star_\Omega\Omega$ & $-\star_\Omega d^{27}_{14}\gamma$ &\\[7pt]
	\end{tabular}
	\caption{Exterior derivative formul{\ae}}\label{extderfor}
\end{table}
That such a table must exist follows from the torsion--freeness which is equivalent to finding coordinates $x_1,\ldots,x_7$ with $\partial_{x_i}\Omega(x)=0$ around any $x\in M$. Put differently, the $\Gt$--structure $(M,\Omega)$ locally osculates to first order to the flat structure on $\R^7$. Therefore, the exterior derivative of expressions such as $\star_\Omega(\alpha\wedge\Omega)$ only depends on the $1$--jet of the reference forms involved, e.g.\ $\alpha$. The operators $d^p_q$ are obtained by compounding $d$ with suitable $\Gt$--equivariant maps. As an example, $d\alpha=\star_\Omega(\dot\alpha\wedge\star_\Omega\Omega)\oplus\dot\beta$ for $\dot\alpha\in\Omega_7$ and $\dot\beta\in\Omega_{14}$, where $\dot\alpha=\star_\Omega(d\alpha\wedge\star_\Omega\Omega)/3$ and $\dot\beta=[d\alpha]_{14}$. Working out the identities of Table~\ref{extderfor} is then a matter of computation using the algebraic formul{\ae}~\eqref{morecontr},~\eqref{omegaisom1} and~\eqref{omegaisom2}. Our precise definition of the operators $d^p_q$ can be found in Appendix~\ref{appendix}.

\begin{rem}
The operators $d^p_q$ and $d^q_p$ are formally adjoint to each other, i.e.
\ben
\langle d^p_q\sigma_p,\sigma_q\rangle_{L^2_{\Omega}}=\langle\sigma_p,d^q_p\sigma_q\rangle_{L^2_{\Omega}}
\ee
for any $\sigma_p\in\Omega_p$ and $\sigma_q\in\Omega_q$.
\end{rem}

\begin{ex}
With Table~\ref{extderfor} at hand we compute the (co-)differential of $\dot\Omega=\dot f\Omega\oplus\star_\Omega(\dot\alpha\wedge\Omega)\oplus\dot\gamma\in\Omega^3(M)$ and find
\begin{eqnarray}
d\dot\Omega & = & \tfrac{4}{7}d^7_1\dot\alpha\star_\Omega\Omega\oplus(d^1_7\dot f+\tfrac{1}{2}d^7_7\dot\alpha+\tfrac{1}{4}d^{27}_7\dot\gamma)\wedge\Omega\oplus\star_\Omega(d^7_{27}\dot\alpha+d^{27}_{27}\dot\gamma)\nonumber\\
\delta_\Omega\dot\Omega & = & \star_\Omega\big((-d^1_7\dot f-\tfrac{2}{3}d^7_7\dot\alpha+\tfrac{1}{3}d^{27}_7\dot\gamma)\wedge\star_\Omega\Omega\big)\oplus d^7_{14}\dot\alpha+d^{27}_{14}\dot\gamma.\label{deltapq}
\end{eqnarray}
\end{ex}

Using Table~\ref{extderfor}, $d^2=0$ implies the following second--order identities of Table~\ref{secordid}.

\medskip

\begin{table}[htb]
	\begin{tabular}{l@{\quad\quad}l@{\quad\quad}l@{\quad\quad}l}
	& $d^7_7d^1_7=0$ & $d^7_{14}d^1_7=0$ &\\[2pt]
	$d^7_1d^7_7=0$ & $d^{14}_7d^7_{14}=\tfrac{2}{3}d^7_7d^7_7$ & $d^7_{14}d^7_7+2d^{27}_{14}d^7_{27}=0$ & $3d^{14}_{27}d^7_{14}+d^7_{27}d^7_7=0$\\[2pt]
	& $d^{27}_7d^7_{27}=d^7_7d^7_7+\tfrac{12}{7}d^1_7d^7_1$ & & $2d^{27}_{27}d^7_{27}-d^7_{27}d^7_7=0$\\[2pt]
	$d^7_1d^{14}_7=0$ & $d^7_7d^{14}_7+2d^{27}_7d^{14}_{27}=0$ & & $d^7_{27}d^{14}_7+4d^{27}_{27}d^{14}_{27}=0$\\[2pt]
	& $3d^{14}_7d^{27}_{14}+d^7_7d^{27}_7=0$ & $d^7_{14}d^{27}_7+4d^{27}_{14}d^{27}_{27}=0$ &\\[2pt]
	& $2d^{27}_7d^{27}_{27}-d^7_7d^{27}_7=0$ & &\\[2pt]
	\end{tabular}
	\caption{Second order identities}\label{secordid}
\end{table}

We will also need the Laplacians $\Delta_\Omega\sigma_p$, $\sigma_p\in\Omega_p$. These are given in Table~\ref{laplace}.

\medskip

\begin{table}[htb]
	\begin{tabular}{lcl}
	\hspace{5cm}$\Delta_\Omega f$ & = & $d^7_1d^1_7f$\\[2pt]
	\hspace{5cm}$\Delta_\Omega\alpha$ & = & $\big(d^7_7d^7_7+d^1_7d^7_1\big)\alpha$\\[2pt]	
	\hspace{5cm}$\Delta_\Omega\beta$ & = & $(\tfrac{5}{4}d^7_{14}d^{14}_7+d^{27}_{14}d^{14}_{27})\beta$\\[2pt]	
	\hspace{5cm}$\Delta_\Omega\gamma$ & = & $\big(\tfrac{7}{12}d^7_{27}d^{27}_7+d^{14}_{27}d^{27}_{14}+(d^{27}_{27})^2\big)\gamma$\\[2pt]
	\end{tabular}
	\caption{Laplacians}\label{laplace}
\end{table}

\begin{ex}
If $\dot\Omega=\dot f\Omega\oplus\star_\Omega(\dot\alpha\wedge\Omega)\oplus\dot\gamma$, then
\ben
\Delta_\Omega\dot\Omega=\Delta_\Omega\dot f\cdot\Omega\oplus\star_\Omega(\Delta_\Omega\dot\alpha\wedge\Omega)\oplus\Delta_\Omega\dot\gamma.
\ee
\end{ex}

\medskip

We finish this section with some material on $\SU(3)$--representation theory. This will become relevant later, when we compute the symbols of various differential operators on a manifold equipped with some fixed $\Gt$--structure. The finer $\SU(3)$--representation theory will enable us to derive properties of these symbols such as invertibility and definiteness.

\bigskip

We pick a unit vector $\xi \in \Lambda^1$ in the vector representation of $\Gt$. Since the unit sphere $S^6$ is diffeomorphic with $\Gt/\SU(3)$, $\xi$ gives rise to an $\SU(3)$--representation over $\xi^{\perp}$, namely the real representation underlying the complex vector representation $\C^3$. In particular, $\xi^\perp$ carries a complex structure. In terms of forms, the group $\SU(3)$ can be regarded as the stabiliser of a non--degenerate $2$--form $\omega\in\Lambda^2\xi^\perp$ and a complex volume form $\Psi=\psi_++i\psi_-\in\Lambda^{3,0}\xi^{\perp}$. These forms relate to $\Omega$ and $\star_\Omega\Omega$ via
\begin{eqnarray}
\Omega & = & \omega\wedge\xi+\psi_+,\label{omdecomp}\\
\star_\Omega\Omega & = & \psi_-\wedge\xi+\frac{1}{2}\omega^2\nonumber.
\end{eqnarray}
In terms of a $\Gt$--frame with $\xi=e_7$ we find $\omega=e^{12}+e^{34}+e^{56}$, $\psi_+=e^{135}-e^{146}-e^{236}-e^{245}$ and $\psi_-=e^{136}+e^{145}+e^{235}-e^{246}$. They satisfy the algebraic relations
\be\label{algrelsu3}
\omega\wedge\psi_\pm=0,\quad\psi_+\wedge\psi_-=\tfrac{2}{3}\omega^3.
\ee 
The decomposition of the exterior algebra over $\xi^\perp$ into irreducibles is given by
\be\label{su3moddecomp}
\lambda^1=\xi^\perp,\quad
\lambda^2=\lambda^2_1\oplus\lambda^2_6\oplus\lambda^2_8,\quad
\lambda^3=\lambda^3_{1+}\oplus\lambda^3_{1-}\oplus\lambda^3_6\oplus\lambda^3_{12},
\ee
where $\lambda^i:=\Lambda^i \xi^\perp$. As above the numerical subscript keeps track of the dimension. We also use these subscripts to denote the corresponding components of a form, e.g.\ $\gamma\in\lambda^3$ can be decomposed into the direct sum $\gamma=\gamma_{1+}\oplus\gamma_{1-}\oplus\gamma_6\oplus\gamma_{12}$. The two trivial representations $\lambda^3_{1\pm}$ are spanned by $\psi_+$ and $\psi_-$ respectively, while $\lambda^2_8$ corresponds to the Lie algebra of $\mf{su}(3)$ sitting inside $\mf{so}(6)\cong\lambda^2$. More importantly for our purposes we can consider the decomposition of the exterior algebra over $\R^7$ into $\SU(3)$--irreducibles. Here, we shall denote by $\mb{(n)}^p_q$ the $n$--dimensional irreducible $\SU(3)$--representation inside $\Lambda^p_q$. Then
\ben
\begin{array}{l}
\Lambda^1\cong(\mb{1})^1_7\oplus(\mb{6})^1_7,\quad\Lambda^2\cong(\mb{1})^2_7\oplus(\mb{6})^2_7\oplus(\mb{6})^2_{14}\oplus(\mb{8})^2_{14},\\[5pt]\Lambda^3\cong(\mb{1})^3_1\oplus(\mb{1})^3_7\oplus(\mb{6})^3_7\oplus(\mb{1})^3_{27}\oplus(\mb{6})^3_{27}\oplus(\mb{8})^3_{27}\oplus(\mb{12})^3_{27},
\end{array}
\ee
so that no confusion shall occur. The decomposition of $\Lambda^3$ is of particular importance for the sequel. The occuring modules can be characterised as follows: 
\ben
\begin{array}{lcl}
(\mb{1})^3_1 & = & \{a(\omega\wedge\xi+\psi_+)\,|\,a\in\R\},\\[5pt]
(\mb{1})^3_7 & = & \{b\psi_-\,|\,b\in\R\},\\[5pt]
(\mb{1})^3_{27} & = & \{c(-4\omega\wedge\xi+3\psi_+)\,|\,c\in\R\},\\[5pt]
(\mb{6})^3_7 & = & \{(X\llcorner\psi_-)\wedge\xi+(X\llcorner\omega)\wedge\omega\,|\,X\in\xi^\perp\},\\[5pt]
(\mb{6})^3_{27} & = & \{(Y\llcorner\psi_-)\wedge\xi-(Y\llcorner\omega)\wedge\omega\,|\,Y\in\xi^\perp\},\\[5pt]
(\mb{8})^3_{27} & = & \{\beta_8\wedge\xi\,|\,\beta_8\in\lambda^2_8\}.
\end{array}
\ee
For instance, $Y\in\xi^\perp\mapsto A(Y)=(Y\llcorner\psi_-)\wedge\xi-(Y\llcorner\omega)\wedge\omega\in\Lambda^3$ is a linear isomorphism onto its image. Further, $(Y\llcorner\psi_-)\wedge\psi_+=-(Y\llcorner\psi_+)\wedge\psi_-$ so that the algebraic relations~\eqref{algrelsu3} readily imply that $A(Y)\wedge\Omega=0$, $A(Y)\wedge\star_\Omega\Omega=0$, i.e.\ $\im A\subset\Lambda^3_{27}$. Summarising, we can write any $\dot\Omega\in\Lambda^3$ as
\begin{eqnarray}\label{finer_decomp}
\dot\Omega & = & [\dot\Omega]_1\oplus[\dot\Omega]_7\oplus[\dot\Omega]_{27}\nonumber\\
& = &\big[\dot a(\omega\wedge\xi+\psi_+)\big]\oplus\big[\dot b\psi_-+(\dot X\llcorner\psi_-)\wedge\xi+(\dot X\llcorner\omega)\wedge\omega\big]\nonumber\\
&   & \oplus\big[\dot c(-4\omega\wedge\xi+3\psi_+)+(\dot Y\llcorner\psi_-)\wedge\xi-(\dot Y\llcorner\omega)\wedge\omega+\dot\beta_8\wedge\xi+\dot\gamma_{12}\big]\label{su3decomp}
\end{eqnarray}
for constants $\dot a,\,\dot b,\,\dot c\in\R$, vectors $\dot X,\,\dot Y\in\xi^\perp$ and forms $\dot\beta_8\in\lambda^2_8$, $\dot\gamma_{12}\in\lambda^3_{12}$. In particular, decomposing $\dot\Omega=\dot\beta\wedge\xi+\dot\gamma$,  where $\dot\beta$ and $\dot\gamma$ are the uniquely determined $2$-- and $3$--forms in $\Lambda^*\xi^{\perp}$ \st $\xi\llcorner\dot\beta,\dot\gamma=0$, we obtain
\begin{eqnarray}
\dot\beta & = & (\dot a-4\dot c)\omega\oplus(\dot X+\dot Y)\llcorner\psi_-\oplus\dot\beta_8\label{betadecomp}\\
\dot\gamma& = & (\dot a+3\dot c)\psi_+\oplus\dot b\psi_-\oplus\big((\dot X-\dot Y)\llcorner\omega\big)\wedge\omega\oplus\dot\gamma_{12}\label{gammadecomp}.
\end{eqnarray}
Thus $\dot\beta_1=(\dot a-4\dot c)\omega$ etc. For later applications, we need for $X\in\xi^\perp$ the identities
\be\label{7star}
\star_\Omega\big((X\llcorner\psi_-)\wedge\Omega\big) =X\llcorner\psi_-+2X\wedge\xi
\ee
and
\be\label{su3isom}
g_\Omega(X\llcorner\psi_-,X\llcorner\psi_-)=2g_\Omega(X,X).
\ee
We prove~\eqref{7star} along the lines of~\eqref{7contrOm}, while~\eqref{su3isom} uses the transitive and isometric action of $SU(3)$ on $S^5$. Hence, up to a rotation we may assume that $X=|X|e_1$. Similarly, the transitive action of $G_2$ on $S^6$ implies
\be\label{omegaisom1}
\langle\xi\wedge\Omega,\xi\wedge\Omega\rangle=4\langle\xi,\xi\rangle,\quad\langle\xi\wedge\star_\Omega\Omega,\xi\wedge\star_\Omega\Omega\rangle=3\langle\xi,\xi\rangle
\ee
for all $\xi\in\Lambda^1$. Furthermore, in conjunction with~\eqref{star2forms} and~\eqref{morecontr} we note the useful formul{\ae}
\begin{eqnarray}\label{omegaisom2}
\langle\tau^2_7\wedge\star_\Omega\Omega,\tau^2_7\wedge\star_\Omega\Omega\rangle &=& 3\langle\tau^2_7,\tau^2_7\rangle,\nonumber\\
\langle\tau^2_7\wedge\Omega,\tau^2_7\wedge\Omega\rangle &=& 4\langle\tau^2_7,\tau^2_7\rangle,\\
\langle\tau^3_7\wedge\Omega,\tau^3_7\wedge\Omega\rangle &=& 4\langle\tau^3_7,\tau^3_7\rangle\nonumber
\end{eqnarray}
for all $\tau^p_q\in\Lambda^p_q$.
%
%
%
%
%
\section{The Dirichlet energy functional $\mc{D}$}\label{def}

In this section we introduce the Dirichlet energy functional and study some basic properties. In particular, we compute its first variation. 

\begin{dfn}
The {\em Dirichlet energy functional} $\mc{D}:\Omega_+^3(M)\to\R_{\geq 0}$ is defined by
\ben
\mc{D}(\Omega) = \frac{1}{2}\int_M\big(|d\Omega|_{g_\Omega}^2+|d\Theta(\Omega)|_{g_\Omega}^2\big)\vol_\Omega.
\ee
\end{dfn}

\begin{rem}
Using the $L^2$--inner product and integration by parts we may also write
\ben
\mc{D}(\Omega)=\frac{1}{2}(\|d\Omega\|^2_{L^2_\Omega}+\|\delta_\Omega\Omega\|^2_{L^2_\Omega})=\frac{1}{2}\big\langle\Delta_\Omega\Omega,\Omega\big\rangle_{L^2_\Omega}.
\ee
\end{rem}

\begin{prp}
{\rm(i)} The functional $\mc{D}$ is invariant under orientation preserving diffeomorphisms, i.e.\ $\mc{D}(\varphi^*\Omega)=\mc{D}(\Omega)$ for all $\varphi\in\diff(M)_+$, $\Omega\in\Omega_+^3(M)$.

{\rm(ii)} For $\lambda\in\R_{>0}$, $\mc{D}(\lambda\Omega)=\lambda^\frac{5}{3}\mc{D}(\Omega)$, i.e.\ $\mc{D}$ is positively homogeneous.
\end{prp}
\begin{prf}
The first assertion follows directly from $\star_{\varphi^*\Omega}=\varphi^*\star_\Omega\varphi^{-1*}$ for $\varphi\in\diff(M)_+$. Secondly, we recall that $\widetilde\Omega=\lambda\Omega\in\Omega^3_+(M)$ if $\lambda>0$ and $\Omega\in\Omega_+^3(M)$. Now a $\Gt$--frame $\{e_i\}$ for $\Omega$ gives the $\Gt$--frame $\{f_i=\lambda^{-1/3}e_i\}$ for $\widetilde\Omega$. Its dual basis is $\{f^i=\lambda^{1/3}e_i\}$. Hence, $\vol_{\lambda\Omega}=f^1\wedge\ldots\wedge f^7=\lambda^\frac{7}{3}\vol_\Omega$, while for the metric $g_{\widetilde\Omega}$ induced on $\Lambda^p$, we find $g_{\lambda\Omega}=f_{i_1\ldots i_p}\otimes f_{i_1\ldots i_p}=\lambda^{-\frac{2p}{3}}g_\Omega$. Hence $|d\widetilde\Omega|^2_{g_{\widetilde\Omega}}=\lambda^{-\frac{2}{3}}|d\Omega|^2_{g_\Omega}$. To compute $|\delta_{\widetilde\Omega}\widetilde\Omega|^2_{g_{\widetilde\Omega}}=|d\star_{\lambda\Omega}\lambda\Omega|^2_{g_{\lambda\Omega}}$ we observe that considered as an operator $\Omega^p(M)\to\Omega^{7-p}(M)$,
\ben
\star_{\lambda\Omega}\alpha^p=\lambda^\frac{7-2p}{3}\star_\Omega\alpha^p,
\ee
whence $|d\star_{\lambda\Omega}\lambda\Omega|^2_{g_{\lambda\Omega}}=\lambda^{-\frac{2}{3}}|d\star_\Omega\!\Omega|^2_{g_\Omega}$.
\end{prf}

\begin{cor}\label{critker}
The space $\mc{X}$ of critical points of $\mc{D}$ is acted on by $\diff(M)_+$ and is given by
\ben
\mc{X}=\{\Omega\in\Omega^3_+(M)\,|\,d\Omega=0,\,\delta_\Omega\Omega=0\},
\ee
the torsion--free positive $3$--forms on $M$, which are the absolute minima of $\mc{D}$.
\end{cor}
\begin{prf}
The first claim follows from diffeomorphism invariance. Secondly, we can apply Euler's formula for homogeneous functions to get
\be\label{positivity}
D_\Omega\mc{D}(\Omega)=\frac{5}{3}\mc{D}(\Omega)=\frac{5}{6}\langle\Delta_\Omega\Omega,\Omega\rangle_{L^2_\Omega}\geq0.
\ee
Equality holds precisely if $\Delta_\Omega\Omega=0$, i.e.\ $d\Omega=0$ and $\delta_\Omega\Omega=0$. Hence, if $\Omega$ is critical, then in particular $D_\Omega\mc{D}(\Omega)=0$ and therefore $\Omega$ is torsion--free.
\end{prf}

Next we compute the first variation of $\mc{D}$. To that end we introduce the following piece of notation. Let $E$ be some vector bundle and $A:\Omega_+^3(M)\to C^\infty(E)$ a differential operator. We write $\dot A_\Omega$ for the linearisation of $A$ at $\Omega\in\Omega_+^3(M)$ evaluated on some $3$--form $\dot\Omega$ tangent to $\Omega$, i.e.\
\ben
\dot A_\Omega:=D_\Omega A(\dot\Omega).
\ee

We illustrate this convention by two examples which will be needed later.

\begin{ex}
(i) Consider the non--linear, homogeneous map
\ben
\Theta:\Omega_+^3(M)\to\Omega^4(M),\quad\Omega\mapsto\star_\Omega\Omega.
\ee
Further, for a fixed $\Gt$--structure $\Omega\in\Omega^3_+(M)$ we define the linear, self--adjoint isomorphism
\ben
p_\Omega:\Omega^3(M)\to\Omega^3(M),\quad\dot\Omega\mapsto\tfrac{4}{3}[\dot\Omega]_1+[\dot\Omega]_7-[\dot\Omega]_{27}.
\ee
With the concrete $\Gt$--structure in mind we shall simply write $p$. By Prop.~10.3.5 in~\cite{jo00}, 
\be\label{thetadot}
\dot\Theta_\Omega=\star_\Omega p_\Omega\dot{\Omega}.
\ee
In particular, $\dot\Omega=\Omega$ gives $\dot\Theta_\Omega=4\Theta(\Omega)/3$.
 
 \medskip
 
(ii) In continuation of the first example we consider the map $F:\Omega_+^3(M)\to\Omega^3(M)$ defined by $F(\Omega)=\Delta_\Omega\Omega$. Then 
\begin{eqnarray*}
\dot{F}_\Omega & = & \dot{\star}_\Omega d\star_\Omega d\Omega+\star_\Omega d\dot{\star}_\Omega d\Omega+\star_\Omega d\star_\Omega d\dot\Omega-d\dot{\star}_\Omega d\Theta(\Omega)-d\star_\Omega d\dot\Theta(\Omega)\\
& \stackrel{(i)}{=} & \star_\Omega d\star_\Omega d\dot\Omega-d\star_\Omega d\star_\Omega p_\Omega\dot\Omega+\mbox{ terms of lower order in }\dot\Omega\nonumber\\
& = & \delta_\Omega d\dot\Omega+d\delta_\Omega p_\Omega\dot{\Omega}+\mbox{ terms of lower order in }\dot\Omega.
\end{eqnarray*}
\end{ex}

\begin{prp}\label{gradient}
We have
\ben
\dot{\mc{D}}_\Omega=\int_M\dot\Omega\wedge\star_\Omega\big(\delta_\Omega d\Omega+p_\Omega d\delta_\Omega\Omega+q_{\Omega}(\nabla^\Omega\Omega)\big)
\ee
for some quadratic form $q_{\Omega}$ whose coefficients depend smoothly on $\Omega$.
\end{prp}
\begin{prf}
As in the previous example, 
\begin{eqnarray}
\dot{\mc{D}}_\Omega & = & \frac{1}{2}\int_Md\dot\Omega\wedge\star_\Omega d\Omega+d\Omega\wedge\big(\dot\star_\Omega d\Omega+\star_\Omega d\dot\Omega\big)\nonumber\\
& & +\frac{1}{2}\int_Md\dot\Theta_\Omega\wedge\star_\Omega d\Theta(\Omega)+d\Theta(\Omega)\wedge\big(\dot\star_\Omega d\Theta(\Omega)+\star_\Omega d\dot\Theta_\Omega\big)\nonumber\\
& = & \int_M d\dot\Omega\wedge\star_\Omega d\Omega+d\dot\Theta_\Omega\wedge\star_\Omega d\Theta(\Omega)\nonumber\\
& & +\frac{1}{2}\int_Md\Omega\wedge\dot\star_\Omega d\Omega+d\Theta(\Omega)\wedge\dot\star_\Omega d\Theta(\Omega).\label{gradD}
\end{eqnarray}
Now $\Gamma_{\Omega}(d\Omega):\dot\Omega\mapsto\dot\star_\Omega d\Omega$ is a bundle endomorphism of $\Lambda^3T^*M$ depending on $d\Omega$ in a fibrewise linear fashion, so we can consider its fibrewise adjoint $(\Gamma_{\Omega}(d\Omega))^*$. Thus
\ben
\int_Md\Omega\wedge\dot\star_\Omega d\Omega=\langle\Gamma_{\Omega}(d\Omega)(\dot\Omega),\star_\Omega d\Omega\rangle_{L^2_\Omega}=\langle\dot\Omega, (\Gamma_{\Omega}(d\Omega))^*(\star_\Omega d\Omega)\rangle_{L^2_\Omega}.
\ee
We deal with the second term of~\eqref{gradD} in a similar manner. The last line is therefore of the form $\int_M\dot\Omega\wedge q_{\Omega}(\nabla^\Omega\Omega)$ with $q_{\Omega}$ quadratic in the first derivatives of $\Omega$, as asserted. On the other hand, Stokes implies
\begin{eqnarray*}
\int_M d\dot\Omega\wedge\star_\Omega d\Omega+d\dot\Theta_\Omega\wedge\star_\Omega d\Theta(\Omega) & = & \int_M \dot\Omega\wedge d\star_\Omega d\Omega-\dot\Theta_\Omega\wedge d\star_\Omega d\Theta(\Omega)\\
& = &\langle\dot\Omega,\delta_\Omega d\Omega\rangle_{L^2_\Omega}+\langle\star_\Omega\dot\Theta_\Omega,d\delta_\Omega \Omega\rangle_{L^2_\Omega}\\
& = & \langle\dot\Omega,\delta_\Omega d\Omega+p_\Omega(d\delta_\Omega\Omega)\rangle_{L^2_\Omega},
\end{eqnarray*}
whence the assertion.
\end{prf}
%
%
%
%
%
\section{Short--time existence}
Let
\ben
Q:\Omega^3_+(M)\to\Omega^3(M),\quad Q(\Omega)=-\grad\mc{D}(\Omega)
\ee
denote the negative $L^2$--gradient of $\mc{D}$ in the sense of Definition 4.10~\cite{be87}, i.e.\ $\langle Q(\Omega),\dot\Omega\rangle_{L^2_\Omega}=-\dot{\mc{D}}_\Omega$. In view of Proposition~\ref{gradient}, we find
\be\label{qop}
Q(\Omega)=-\delta_\Omega d\Omega-p_\Omega(d\delta_\Omega\Omega)-q_{\Omega}(\nabla\Omega).
\ee

The goal of this section is to prove the existence part of

\begin{thm}\label{exandun}
Given $\Omega_0\in\Omega_+^3(M)$, there exists $\epsilon>0$ and a smooth family $\Omega(t)\in\Omega_+^3(M)$ for $t\in[0,\epsilon]$ \st
\be\label{floweq}
\frac{\partial}{\partial t}\Omega=Q(\Omega),\quad\Omega(0)=\Omega_0.
\ee  
Further, if $\Omega(t)$ and $\Omega'(t)$ are solutions to~\eqref{floweq}, then $\Omega(t)=\Omega'(t)$ whenever defined. Hence $\Omega(t)$ is uniquely defined on a maximal time--interval $[0,T)$ for some $0<T\leq\infty$.
\end{thm}

\begin{rem}
We emphasise that, in contrast to the corresponding results for the Laplacian flow in \cite{brxu11}, \cite{xuye09}, the initial condition $\Omega_0$ is not assumed to be closed.
\end{rem}

\begin{dfn}
We call the negative gradient flow of $\mc{D}$ defined by~\eqref{floweq} the {\em Dirichlet energy flow with initial condition} $\Omega_0\in\Omega_+^3(M)$.
\end{dfn}

\medskip

We will prove short--time existence and uniqueness by invoking the standard theory of quasilinear parabolic equations which we briefly recall, see Chapter 4.4.2~\cite{au97}. Further useful references are Chapter 7.8 in~\cite{lsu67} and Chapter 7.1 in~\cite{ta91}. Consider a Riemannian vector bundle $\big(E,(\cdot\,,\cdot)\big)$. Let $Q_t:C^{\infty}(E)\to C^{\infty}(E)$ be a family of quasilinear, second order differential operators, that is locally, $Q_t(u)(x)\loc\big(a_{\beta}^{\alpha ij}(t,x,u,\nabla u)\partial_i\partial_j u^\beta+b^\alpha(t,x,u,\nabla u)\big)s_\alpha$ for smooth functions $a^{\alpha ij}_\beta$ and $b^\alpha$ and a local basis $\{s_\alpha\}$ of $E$. We say that the induced flow equation
\be\label{standeq}
\frac{\partial}{\partial t}u=Q_t(u),\quad u(0)=u_0
\ee
is {\em strongly parabolic} at $u_0$ if there exists a constant $\lambda>0$ \st the linearisation $D_{u_0}Q_0$ of $Q_0$ at $u_0$ satisfies
\be\label{strongpara}
-\big(\sigma(D_{u_0}Q_0)(x,\xi)v,v\big)\geq\lambda|\xi|^2|v|^2
\ee
for all $(x,\xi)\in TM$, $\xi\not=0$, and $v\in E_x$. Here, the minus sign in~\eqref{strongpara} stems from our definition of the principal symbol. Namely, for a $k$--th order linear differential operator $Q$ we define
\ben
\sigma(Q)(x,\xi)v=\tfrac{i^k}{k!}Q(f^ku)(x)
\ee
for an $f\in C^\infty(M)$ with $f(x)=0$, $d_xf=\xi$ and $u\in C^\infty(E)$ with $u(x)=v$.

\begin{thm}\label{mainanaresult}
If equation~\eqref{standeq} is strongly parabolic at $u_0$, then there exists $\epsilon>0$ and a smooth family $u(t)\in C^{\infty}(E)$, $t\in[0,\epsilon]$ \st 
\ben
\frac{\partial}{\partial t}u=Q_t(u),\quad u(0)=u_0.
\ee
Further, if $u(t)$ and $u'(t)$ are solutions to~\eqref{floweq}, then $u(t)=u'(t)$ whenever defined. 
\end{thm}

\bigskip

Next we investigate the operator $Q$ as given in~\eqref{qop}.

\begin{lem}\label{qanaprop}
The second order non--linear differential operator $Q$ is quasilinear.
\end{lem}
\begin{prf}
For instance, up to composition with the linear map $p$ whose coefficients depend solely on $\Omega$, the second term on the right hand side of~\eqref{qop} can be locally written as
\ben
-d\delta_\Omega\Omega\loc\partial_i\big(\star^{jstuv}_{\Omega,pq}(\Omega)\partial_j(\star^{abc}_{\Omega, stuv}(\Omega)\Omega_{abc})\big)dx^{ipq}.
\ee
Here, $\star^{abc}_{\Omega,stuv}(\Omega)$ denote the coefficients of $\star_\Omega:\Omega^3(M)\mapsto\Omega^4(M)$ \wrt local coordinates $x_1,\ldots,x_7$ etc. Hence $Q$ is linear in its highest (i.e.~second) order derivatives.
\end{prf}

\begin{lem}\label{linQ}
The principal symbol $\sigma(D_\Omega Q)(x,\xi):\Lambda^3T_x^*M\to\Lambda^3T_x^*M$ of the linearisation of $Q$ at $\Omega\in\Omega_+^3(M)$ is given by
\ben
\sigma(D_\Omega Q)(x,\xi)\dot\Omega = -\xi\llcorner(\xi\wedge\dot\Omega)-p_\Omega\big(\xi\wedge\big(\xi\llcorner p_\Omega(\dot\Omega)\big)\big).
\ee
Moreover, the symbol is negative semi--definite.
\end{lem}
\begin{prf}
As the principal symbol involves highest order terms only, we only need to linearise the expression
\ben
-\delta_\Omega d\Omega-p_\Omega(d\delta_\Omega\Omega)=Q(\Omega)-\mbox{ terms of lower order in }\Omega.
\ee
In our convention, $\sigma(d)(x,\xi)\dot\Omega=i\xi\wedge\dot\Omega$ and $\sigma(\delta_\Omega)(x,\xi)\dot\Omega=-i\xi\llcorner\dot\Omega$. Hence, from~\eqref{thetadot} and the standard symbol calculus we get the asserted symbol. Further,
\begin{eqnarray*}
-g_\Omega\big(\sigma(D_\Omega Q)(x,\xi)\dot\Omega,\dot\Omega\big) & = & g_\Omega\big( \xi\llcorner(\xi\wedge\dot\Omega)+p_\Omega\big(\xi\wedge\big(\xi\llcorner p_\Omega(\dot\Omega)\big)\big)
,\dot\Omega\big)\\
& = & |\xi\wedge\dot\Omega|^2_\Omega+|\xi\llcorner p_\Omega(\dot\Omega)|^2_\Omega\\
& \geq & 0
\end{eqnarray*}
so that $\sigma(D_\Omega Q)(x,\xi)$ is negative semi--definite.
\end{prf}

\begin{rem}
For $\varphi\in\diff(M)_+$ we have
\be\label{commutQ}
\varphi^*Q(\Omega)=Q(\varphi^*\Omega)
\ee
since $\mc{D}\circ\varphi^*=\mc{D}$ and $\langle\dot\Omega_0,\dot\Omega_1\rangle_{L^2_{\varphi^*\Omega}}=\langle\varphi^{-1*}\dot\Omega_0,\varphi^{-1*}\dot\Omega_1\rangle_{L^2_\Omega}$. Because of this diffeomorphism invariance we cannot expect the principal symbol to be negative definite. Indeed, using the decomposition $\dot\Omega=\dot\beta\wedge\xi+\dot\gamma$ as in Section~\ref{reptheory}, $g_\Omega\big(\sigma(D_\Omega Q)(x,\xi)\dot\Omega,\dot\Omega\big)=0$  implies $\dot\gamma=0$. From~\eqref{gammadecomp} we deduce $\dot a=-3\dot c$, $\dot b=0$, $\dot X_0=\dot Y_0$ and $\dot\gamma_{12}=0$. Consequently, one then gets that $\dot\beta_8=-\xi\llcorner p(\dot\Omega)=0$, whence 
\ben
\ker\sigma(D_\Omega Q)(x,\xi)=\{(\dot v\omega+\dot V\llcorner\psi_-)\wedge\xi\,|\,\dot v\in\R,\,\dot V\in\xi^\perp\}
\ee
by~\eqref{su3decomp}.
\end{rem}

\medskip

In order to apply Theorem~\ref{mainanaresult} we use so--called DeTurck's trick which was orginally invented for Ricci flow~\cite{dt83}. Given a family of diffeomorphisms $\partial_t\varphi_t=X_t\circ\varphi_t$ induced by a (time--dependent) vector field $X_t$ on $M$, differentiating~\eqref{commutQ} yields the intertwining formula
\be\label{infcommutQ}
\mc{L}_X\big(Q(\Omega)\big)=D_\Omega Q(\mc{L}_X\Omega).
\ee
Here $\mc{L}_X$ denotes Lie derivative \wrt $X$. While the left hand side of~\eqref{infcommutQ} is of first order in $X$, the right hand side is of third order. Passing to symbol level implies
\be\label{imageinkernel}
\sigma(D_\Omega Q)(x,\xi)\circ\sigma(X\mapsto\mc{L}_X\Omega)(x,\xi)=0.
\ee
In this way, we can conceive the symbol of the map
\be\label{lambda}
\Omega\in\Omega_+^3(M)\mapsto X(\Omega)\in C^{\infty}(TM)\mapsto\Lambda(\Omega)=\mc{L}_{X(\Omega)}\Omega\in\Omega^3(M)
\ee
(where the vector field $X(\Omega)$ depends non--trivially on the $1$--jet of $\Omega$) as a kind of projection to the kernel of $\sigma(D_\Omega Q)$. One therefore expects the symbol of the modified operator
\be\label{pertop}
\widetilde{Q}(\Omega)=Q(\Omega)+\Lambda(\Omega)
\ee
to have trivial kernel for a suitably chosen vector field. For a fixed $\bar\Omega\in\Omega_+^3(M)$ we take
\be\label{DTvf}
X_{\bar\Omega}:\Omega^3(M)\to\Omega^1(M),\quad X_{\bar\Omega}(\dot\Omega):=-(\delta_{\bar\Omega}\dot\Omega)\llcorner\bar\Omega,
\ee
where we contract and dualise \wrt the metric $g_{\bar\Omega}$. We think of $X_{\bar\Omega}$ as a first order, linear differential operator. Subsequently we write $\Lambda_{\bar\Omega}$ and $\widetilde Q_{\bar\Omega}$ in~\eqref{pertop} to emphasise the choice of $\bar\Omega$.

\bigskip

To give some motivation concerning the definition of $X_{\bar\Omega}$ we introduce the operator
\be\label{lambdastern}
\lambda^*_{\bar\Omega}:C^\infty(TM)\to\Omega^3(M),\quad X\mapsto\mc{L}_X\bar\Omega.
\ee 
We consider $\lambda^*_{\bar\Omega}$ to be the formal $L^2$--adjoint \wrt $\langle\cdot\,,\cdot\rangle_{L^2_{\bar\Omega}}$ of $\lambda_{\bar\Omega}$, i.e.
\begin{eqnarray*}
\langle\lambda^*_{\bar\Omega}(X),\dot\Omega\rangle_{L^2_{\bar\Omega}} & = & \langle X\llcorner d\bar\Omega+d(X\llcorner\bar\Omega),\dot\Omega\rangle_{L^2_{\bar\Omega}}\\
& = & -\langle d\bar\Omega,\dot\Omega\wedge X\rangle_{L^2_{\bar\Omega}}+\langle\bar\Omega,\delta_{\bar\Omega}\dot\Omega\wedge X\rangle_{L^2_{\bar\Omega}}\\
& = & \langle-\dot\Omega\llcorner d\bar\Omega+\delta_{\bar\Omega}\dot\Omega\llcorner \bar\Omega,X\rangle_{L^2_{\bar\Omega}},
\end{eqnarray*}
whence
\be\label{lambdaomega}
\lambda_{\bar\Omega}(\dot\Omega)=-X_{\bar\Omega}(\dot\Omega)-\dot\Omega\llcorner d\bar\Omega=-X_{\bar\Omega}(\dot\Omega)+\mbox{ terms of lower order in }\dot\Omega.
\ee
In analogy with the decomposition of symmetric $2$--tensors into a divergence free part and a part tangential to the $\diff(M)_+$--orbit of some given Riemannian metric we have:

\begin{prp}\label{divsplit}
For any $\dot\Omega\in\Omega^3(M)$ there exist $X\in C^\infty(TM)$ and $\dot\Omega_0\in\Omega^3(M)$ with $\lambda_{\bar\Omega}(\dot\Omega_0)=0$ \st we have an $L^2_{\bar\Omega}$--orthogonal decomposition 
\be\label{divcomp}
\dot\Omega=\dot\Omega_0\oplus\mc{L}_X\bar\Omega.
\ee
\end{prp}
\begin{prf}
Put $L=\lambda_{\bar\Omega}\lambda^*_{\bar\Omega}:C^\infty(TM)\to C^\infty(TM)$. If we can solve $\lambda_{\bar\Omega}(\dot\Omega)=L(X)$ for some $X\in C^{\infty}(TM)$, then taking $\dot\Omega_0=\dot\Omega-\mc{L}_X\bar\Omega$ yields the desired splitting. Since $L$ is symmetric, such an $X$ exists \iff $\lambda_{\bar\Omega}(\dot\Omega)\in(\ker L)^\perp$. But $L(Y)=0$ implies $\lambda^*_{\bar\Omega}(Y)=0$, whence $\langle\lambda_{\bar\Omega}(\dot\Omega),Y\rangle_{L^2_{\bar \Omega}}=0$ for all $Y\in\ker L$.
\end{prf}

\begin{rem}
(i) The condition $\lambda_{\bar \Omega}(\dot \Omega)=0$ should be viewed as a gauge--fixing condition, i.e.~a choice of a local slice to the $\diff(M)_+$--action near $\bar \Omega$. If $\bar \Omega$ is torsion--free, then $\lambda_{\bar \Omega}(\dot \Omega) = X_{\bar \Omega}(\dot \Omega) = (\delta_{\bar \Omega}\dot\Omega) \llcorner \bar \Omega$, hence $\lambda_{\bar \Omega}(\dot \Omega) = 0$ \iff $[\delta_{\bar \Omega}\dot\Omega]_7=0$. This is precisely the gauge--fixing condition considered by Joyce in \cite{jo00}, see also the remark following Corollary~\ref{pert-unpert}.

\medskip

(ii) The vector field $X$ in the decomposition~\eqref{divcomp} is unique if there are no non--trivial infinitesimal automorphisms of $\bar\Omega$, i.e.\ vector fields $X$ \st $\mc{L}_X\bar\Omega=0$. This holds for holonomy $\Gt$--manifolds as Ricci--flatness implies any Killing field $X$ to be parallel, so that the holonomy is contained in $\SU(3)$ unless $X=0$. Note further that a generic $\Gt$--form has no infinitesimal automorphisms as these are automatically Killing for $g_{\bar\Omega}$ and a generic metric has no Killing vector fields~\cite{eb68}. An example for a non--generic $3$--form is provided by the direct product $M=S^1\times CY^6$ of $S^1$ with coordinate vector field $X=\partial_t$ and an almost Calabi--Yau manifold $(CY^6,\omega,\psi_+)$ (i.e.\ $\omega$ is a non--degenerate $2$--form and $\psi_+$ is a $3$--form of special algebraic type). By~\eqref{omdecomp}, $\bar\Omega=dt\wedge\omega+\psi_+$ is a $\Gt$--form and $\mc{L}_X\bar\Omega=0$.  
\end{rem}

\medskip

Coming back to the mainstream development, we establish strong parabolicity for the flow equation
\be\label{ddflow}
\frac{\partial}{\partial t}\widetilde\Omega=\widetilde Q_{\bar\Omega}(\widetilde\Omega),\quad\widetilde\Omega(0)=\Omega_0.
\ee

\begin{lem}\label{qtildeoperator}
The equation~\eqref{ddflow} is strongly parabolic for $\Omega_0$ sufficiently $C^2$--close to $\bar\Omega$.
\end{lem}
\begin{prf}
Since $X_{\bar\Omega}$ is linear in $\Omega$, we find for the linearisation
\ben
\dot\Lambda_{\bar\Omega}=d\big(X_{\bar\Omega}(\dot\Omega)\llcorner\bar\Omega\big)+\mbox{ lower order terms in }\dot\Omega
\ee
by virtue of Cartan's formula, whence
\begin{eqnarray*}
\sigma(D_{\bar\Omega}\Lambda)(x,\xi)\dot\Omega & = & i\xi\wedge\big(\sigma(X_{\bar\Omega})(x,\xi)\dot\Omega\llcorner\bar\Omega\big)\\
& = & -\xi\wedge\big((\xi\llcorner\dot\Omega)\llcorner\bar\Omega\llcorner\bar\Omega\big).
\end{eqnarray*}
Assume without loss of generality that $|\xi|_{\bar\Omega}=1$. Decomposing $\dot\Omega=\beta\wedge\xi+\gamma$ as above we deduce from~\eqref{7contrOm}
\ben
\sigma(D_{\bar\Omega}\Lambda)(x,\xi)\dot\Omega=-\xi\wedge\big((\beta\llcorner\bar\Omega)\llcorner\bar\Omega\big)=-3\xi\wedge[\beta]_7.
\ee
Bearing~\eqref{star2forms} in mind, the projection of $\beta$ onto $\Lambda^2_7$ is given by
\begin{eqnarray*}
[\beta]_7 & = & \tfrac{1}{3}\big(\beta+\star_{\bar\Omega}(\beta\wedge\bar\Omega)\big)\\
& = & \beta_1\oplus\tfrac{1}{3}\big(\beta_6+\star_{\bar\Omega}(\beta_6\wedge\bar\Omega)\big)\\
& \stackrel{\eqref{betadecomp},\,\eqref{7star}}{=} & \beta_1\oplus\tfrac{2}{3}\big((X_0+Y_0)\llcorner\psi_-+(X_0+Y_0)\wedge\xi\big).
\end{eqnarray*}
Consequently $\xi\wedge[\beta]_7=\xi\wedge\big(\beta_1+\frac{2}{3}(X_0+Y_0)\llcorner\psi_-\big)$ so that using 
$$
|\xi \wedge \gamma|_{\bar\Omega}^2 = |\gamma|_{\bar \Omega}^2 \quad \text{and} \quad g_{\bar \Omega} (\xi \wedge [\beta]_7, \xi \wedge \beta + \gamma) = | [\beta]_7|_{\bar \Omega}^2
$$
the computation from Lemma~\ref{linQ} implies
\be\label{strongellip}
g_{\bar\Omega}\big(\sigma(D_{\bar\Omega}\widetilde Q)(x,\xi)\dot\Omega,\dot\Omega\big)=-|\gamma|_{\bar \Omega}^2-|\xi\llcorner p_{\bar\Omega}(\dot\Omega)|_{\bar \Omega}^2-3|\beta_1|_{\bar\Omega}^2-2|(X_0+Y_0)\llcorner\psi_-|_{\bar \Omega}^2.
\ee
Now $\xi\llcorner p_{\bar\Omega}(\dot\Omega)=\sigma\oplus-\beta_8$ with $g_{\bar\Omega}(\sigma,\beta_8)=0$ in view of the decomposition in \eqref{finer_decomp}, while by~\eqref{su3moddecomp}, $|\beta|_{\bar\Omega}^2=|\beta_1|_{\bar\Omega}^2+|\beta_6|_{\bar\Omega}^2+|\beta_8|_{\bar\Omega}^2$. But~\eqref{betadecomp} gives $|\beta_6|_{\bar\Omega}^2=|(X_0+Y_0)\llcorner\psi_-|_{\bar\Omega}^2$, whence
\ben
-g_{\bar\Omega}\big(\sigma(D_{\bar\Omega}\widetilde Q)(x,\xi)\dot\Omega,\dot\Omega\big)\geq(|\beta|_{\bar\Omega}^2+|\gamma|_{\bar\Omega}^2)=|\dot\Omega|_{\bar\Omega}^2
\ee
by~\eqref{strongellip}.
\end{prf}

\begin{dfn}
We call the flow associated with~\eqref{ddflow}, the {\em Dirichlet--DeTurck flow at $\bar\Omega$ with initial condition} $\Omega_0$. If the $\Gt$--form $\bar\Omega$ we use to perturb the Dirichlet energy flow is understood we simply speak of the Dirichlet--DeTurck flow with initial condition $\Omega_0$.
\end{dfn}

For $\Omega_0\in\Omega^3_+(M)$ consider the Dirichlet--DeTurck flow at $\Omega_0$ with initial condition $\Omega_0$. By Theorem~\ref{mainanaresult}, the flow $\widetilde\Omega(t)$ exists on some time interval $[0,\epsilon]$. Let $\varphi_t$ be the family of diffeomorphisms determined by
\be\label{forwardflow}
\partial_t\varphi_t=-X_{\Omega_0}\big(\widetilde\Omega(t)\big)\circ\varphi_t,\quad\varphi_0=\Id_M.
\ee
Then $\Omega(t)=\varphi_t^*\widetilde\Omega(t)$ is a solution to the Dirichlet energy flow~\eqref{floweq} with same initial condition $\Omega_0$ for
\begin{eqnarray*}
\frac{\partial}{\partial t}\Omega & = & \varphi^*_t\big(\frac{\partial}{\partial t}\widetilde\Omega+\mc{L}_{-X_{\Omega_0}(\widetilde\Omega)}\widetilde\Omega\big)\\
& \stackrel{\eqref{pertop}}{=} & \varphi^*_tQ(\widetilde\Omega)\\
& \stackrel{\eqref{commutQ}}{=} & Q(\Omega).
\end{eqnarray*}
Moreover, the initial condition is satisfied as $\Omega(0)=\Id_M^*\Omega_0=\Omega_0$.

\begin{cor}[Existence]\label{pert-unpert}
For any $\Omega_0\in\Omega^3_+(M)$ there exists an $\epsilon>0$ \st the Dirichlet energy flow~\eqref{floweq} exists for $t\in[0,\epsilon]$.
\end{cor}

\begin{rem}
(i) The idea of DeTurck's trick is to break the diffeomorphism invariance by modifying the flow along the $\diff(M)_+$--orbits via the additional term $\Lambda_{\Omega_0}$. To see this happening in a geometrical way, assume for simplicity that $\Omega_0$ is a torsion--free $\Gt$--form. By~\eqref{7contrOm}, $\dot\Lambda_{\Omega_0}=d\big(X_{\Omega_0}(\dot\Omega)\llcorner\Omega_0)=-3d[\delta_{\Omega_0}\dot\Omega]_7$. In particular, $\dot\Lambda_{\Omega_0}=0$ \iff $[\delta_{\Omega_0}\dot\Omega]_7=0$, for 
$$
\langle d[\delta_{\Omega_0}\dot\Omega]_7,\dot\Omega\rangle_{L^2_{\Omega_0}}=\|[\delta_{\Omega_0}\dot\Omega]_7\|^2_{L^2_{\Omega_0}}. 
$$
On the other hand, the tangent space at $\Omega_0$ of the $\diff(M)_+$--orbit $\mc{O}_{\Omega_0}$ is given by $\{\mc{L}_X\Omega_0\,|\,X\in C^\infty(TM)\}$ (cf.\ also Lemma~\ref{orbit} below). Since 
$$
\langle\dot\Omega-\Omega_0,\mc{L}_X\Omega_0\rangle_{L^2_{\Omega_0}}=\langle[\delta_{\Omega_0}\dot\Omega]_7,X\llcorner\Omega_0\rangle_{L^2_{\Omega_0}},
$$ 
the form $\dot\Omega-\Omega_0$ is perpendicular to $T_{\Omega_0}\mc{O}_{\Omega_0}$ \iff $\dot\Lambda_{\Omega_0}=0$. 

\medskip

(ii) To become strongly parabolic after perturbing with $\Lambda_{\Omega_0}$ is a particular feature of the Dirichlet energy flow. In contrast consider the gradient flow of the Hitchin functional $\mc{H}$ restricted to the cohomology class $[\Omega_0]$, cf.\ \eqref{functional}. Upon rescaling $\mc{H}$, the resulting flow is 
\ben
\frac{\partial}{\partial t}\alpha=\delta_{\Omega_0+d\alpha}(\Omega_0+d\alpha),\quad\alpha(0)=0
\ee
for $\alpha$ in a suitably small open neighbourhood of $0\in\Omega^2(M)$ so that $\Omega_0+d\alpha\in\Omega_+^3(M)$. The solutions are in 1--1 correspondence with solutions to the Laplacian flow 
\ben
\frac{\partial}{\partial t}\Omega=F(\Omega)=\Delta_\Omega\Omega,\quad\Omega(0)=\Omega_0.
\ee
Now $\sigma(D_\Omega F)(x,\xi)\dot\Omega=|\xi|^2_\Omega\dot\Omega + \xi\wedge\big(\xi\llcorner([\dot\Omega]_1/3-2[\dot\Omega]_{27})\big)$ and one can easily compute that the kernel is also given by 
$$
K=\ker\sigma(D_\Omega F)(x,\xi)=\{(\dot v\omega+\dot V\llcorner\psi_-)\wedge\xi\,|\,\dot v\in\R,\,\dot V\in\xi^\perp\}. 
$$
Since by~\eqref{imageinkernel}, the symbol of $X\mapsto\mc{L}_X\Omega$ takes values in $K$, DeTurck's trick cannot modify the component $\mr{pr}_{K^\perp}\circ\sigma(D_\Omega F)_{|K^\perp}:K^\perp\to K^\perp$. However, the eigenvectors $\dot\Omega_1=\dot\beta_8\wedge\xi$ and $\dot\Omega_2=\psi_-$ in $K^\perp$ give $g_\Omega\big(\dot\Omega_1,\sigma(D_\Omega F)(x,\xi)\dot\Omega_1\big)=-|\xi|^2_\Omega|\dot\beta_8|^2_\Omega<0$ and $g_\Omega\big(\dot\Omega_2,\sigma(D_\Omega F)(x,\xi)\dot\Omega_2\big)=4|\xi|^2_\Omega>0$ respectively. Hence, the linearisation of $\widetilde F=F+\Lambda$ will be indefinite no matter how the vector field $X$ in~\eqref{lambda} is chosen (though the linearisation of $\widetilde F$ might have trivial kernel). We therefore deal with a heat equation of mixed forwards/backwards type for which short--time existence is in general not expected unless further conditions are imposed. For the Laplacian flow this has recently been achieved in \cite {brxu11} and \cite{xuye09} for closed initial conditions.
\end{rem}
%
%
%
%
%
\section{Uniqueness}
We now settle the uniqueness part of Theorem~\ref{exandun} along the lines of the uniqueness proof for Ricci flow. 

\bigskip

As shown by Corollary~\ref{pert-unpert}, a solution to the Dirichlet--DeTurck flow $\widetilde\Omega(t)$ with initial condition $\Omega_0$ yields a solution to the Dirichlet energy flow $\Omega(t)=\varphi^*_t\widetilde\Omega(t)$ with same initial condition by integrating the time--dependent vector field in~\eqref{forwardflow}. Conversely, substituting $\widetilde\Omega(t)$ by $\varphi^{-1*}_t\Omega(t)$ turns the ordinary differential equation~\eqref{forwardflow} into the partial differential equation
\be\label{backwardflow}
\frac{\partial}{\partial t}\varphi_t=-X_{\Omega_0}\big(\varphi^{-1*}_t\Omega(t)\big)\circ\varphi_t,\quad\varphi_0=\Id_M.
\ee
A curve $\varphi_t\in\diff(M)_+$ which solves~\eqref{backwardflow} for a Dirichlet energy flow solution $\Omega(t)$ with initial condition $\Omega_0$
yields the Dirichlet--DeTurck flow solution $\widetilde\Omega(t)=\varphi_t^{-1*}\Omega(t)$ with same initial condition. Indeed, let $Y_t$ be the time--dependent vector field defined by $Y_t\circ\varphi_t^{-1}=\partial_t\varphi_t^{-1}$. Then differentiating the constant curve $\varphi^{-1}_t\circ\varphi_t(x)=x$ gives
\be\label{vfisgood}
Y_t(x)=-d_{\varphi_t(x)}\varphi^{-1}_t\big(-X_{\Omega_0}\big(\varphi^{-1*}_t\Omega(t)\big)\circ\varphi_t(x)\big)=\varphi^{-1}_{t*}X_{\Omega_0}\big(\varphi^{-1*}_t\Omega(t)\big)(x),
\ee
where for $\varphi\in\diff(M)$ and $X\in C^{\infty}(TM)$,
\ben
(\varphi_*X)(x):=d_{\varphi^{-1}(x)}\varphi\big(X(\varphi^{-1}(x)\big).
\ee
As a consequence, we get 
\begin{eqnarray*}
\frac{\partial}{\partial t}\varphi^{-1*}_t\Omega(t) & = &\varphi^{-1*}_t\big(\frac{\partial}{\partial t}\Omega(t)+\mc{L}_{Y_t}\Omega(t)\big)\\
& = & Q\big(\varphi^{-1*}_t\Omega(t)\big)+\mc{L}_{\varphi_{t*}Y_t}\varphi^{-1*}_t\Omega(t)\\
& = & Q\big(\varphi^{-1*}_t\Omega(t)\big)+\mc{L}_{X_{\Omega_0}(\varphi^{-1*}_t\Omega(t))}\varphi_t^{-1*}\Omega(t)
\end{eqnarray*}
by~\eqref{vfisgood}. We can then deduce uniqueness of the Dirirchlet energy flow from uniqueness of the Dirichlet--DeTurck flow, see Corollary~\ref{DF-unique}.

\begin{rem}
Equation~\eqref{backwardflow} should be considered as an analogue of the harmonic map heat flow
\ben
\frac{\partial}{\partial t}\varphi_t=\tau_{g(t),g_0}(\varphi_t)
\ee
introduced by Eells and Sampson~\cite{eesa64}, albeit with a time--dependent tension field $\tau_{g(t),g_0}(\varphi_t)$. We can think of $\tau_{g(t),g_0}(\varphi_t)$ as a differential operator defined by Riemannian metrics $g(t)$ and $g_0$ on $M$, taking a smooth map $\varphi:\big(M,g(t)\big)\to (M,g_0)$ to a section $\tau_{g(t),g_0}(\varphi)\in C^{\infty}(\varphi^*TM)$. 
\end{rem}

We need to prove short--time existence of a solution
to~\eqref{backwardflow}. By and large, we proceed as in the harmonic map
heat flow case, cf.~\cite{eesa64}. Let
\be\label{finalflow}
P_t=P_{\Omega(t),\Omega_0}:\varphi\in\diff(M)_+\subset C^{\infty}(M,M)\mapsto-\varphi_*X_{\varphi^*\Omega_0}\big(\Omega(t)\big)\circ\varphi\in C^{\infty}(\varphi^*TM).
\ee

Since $\diff(M)_+$ is open in $C^{\infty}(M,M)$, a solution to the flow equation
\ben
\frac{\partial}{\partial t}\varphi_t=P_t(\varphi_t),\quad\varphi_0=\Id_M
\ee
yields the desired solution to~\eqref{backwardflow}.

\medskip

To get formally in a situation to apply Theorem~\ref{mainanaresult}, we first choose an embedding $f_0:M\to\R^n$ and identify $M$ with its image under $f_0$. In particular, all tensors on $M$ pushed forward to $f_0(M)$ will be denoted by the same symbol. Let $\mc{N}\subset\R^n$ be a tubular neighbourhood of $M$ which we think of as an open neighbourhood inside the normal bundle $\pi:\nu M\to M$, the normal bundle taken with respect to the Euclidean metric on $\R^n$. By choosing a fibre metric $h$ and a compatible connection $\nabla^{\nu M}$ on $\nu M$ we obtain the induced metric $\pi^*g_{\Omega_0}+h$ on $\mc{N}$ which we extend to $\R^n$ using a partition of unity. In particular, this makes $f_0$ an isometry. Similarly, we extend $\Omega_0$ by $\pi^*\Omega_0$ to $\mc{N}$ and subsequently to $\R^n$. In this way, the restriction $f^*\Omega_0$ for $f$ in a suitably small open neighbourhood $\mc{U}\subset C^{\infty}(M,\R^n)$ of embeddings close to $f_0$ is still a positive $3$--form on $M$. Consequently, we can extend $P_t$ to an operator
\ben
P_t:\mc{U}\subset C^{\infty}(M,\R^n)\to C^{\infty}(M,\R^n),\quad f\mapsto-df\big(X_{f^*\Omega_0}\big(\Omega(t)\big)\big).
\ee

\begin{lem}
The operator $P_t$ is a quasilinear, second order differential operator.
\end{lem}
\begin{prf}
Let $e_1,\ldots,e_n$ be the standard basis of $\R^n$ and $x^1,\ldots,x^7$ be
local coordinates on $U \subset M$. The components $\star^{ijk}_{opqr}$ of $\star_{f^*\Omega_0}:\Omega^3(U)\mapsto\Omega^4(U)$ depend on the components of $f^*\Omega_0$ given by $\Omega_{0,\alpha\beta\gamma}\partial_{x_l}f^\alpha\partial_{x_m}f^\beta\partial_{x_n}f^\gamma$. Schematically,
\ben
\star_{f^*\Omega_0}\Omega(t)\loc\star^{ijk}_{opqr}(x,\partial_{x_l}f^\alpha)\Omega(t)_{ijk}dx^{opqr}
\ee
so that by the chain rule
\ben
d\star_{f^*\Omega_0}\Omega(t)\loc\big(a^{ij}_{opqrs,\beta}(t,x,\partial_{x_l}f^\alpha)\partial_{x_i}\partial_{x_j}f^\beta+b_{opqrs}(t,x,\partial_{x_l}f^\alpha)\big)dx^{opqrs}
\ee
for smooth coefficients $a^{\ldots}_{\ldots}$ and $b_{\ldots}$. Applying once more $\star_{f^*\Omega_0}$ and contracting the result with $f^*\Omega_0$ leads to
\ben
\big(\star_{f^*\Omega_0}d\star_{f^*\Omega_0}\Omega(t)\big)\llcorner f^*\Omega_0\loc\big(\widetilde a^{ij}_{k,\beta}(t,x,\partial_{x_l}f^\alpha)\partial_i\partial_jf^\beta+\widetilde b_{kl}(t,x,\partial_{x_l}f^\alpha)\big)dx^k.
\ee
Finally, dualising and contracting with $df=\partial_{x_i}f^\gamma dx^i\otimes e_\gamma$ shows that $P_t$ is a quasilinear, second order differential operator. 
\end{prf}

\begin{lem}
There exists $\epsilon>0$ and a smooth family of embeddings $f(t)\in C^{\infty}(M,\R^n)$, $t\in[0,\epsilon]$ \st
\be\label{unicflow}
\frac{\partial}{\partial t}f(t)=P_t\big(f(t)\big),\quad f(0)=f_0.
\ee
Furthermore, $f(t)(M)\subset f_0(M)$ for all $t$.
\end{lem}
\begin{prf}
As above let $\mathcal U \subset C^\infty(M,\R^n)$ be an open neighbourhood of $f_0$. By transversality of $f_0(M)$ to the fibres of the normal bundle we may assume that $f$ is an embedding for $f \in \mathcal{U}$ with $\varphi_f:=\pi\circ f\in\diff(M)_+$, shrinking $\mc{U}$ if necessary. We put $\sigma_f(x):=f(x)-\varphi_f(x) \in \nu_{\varphi_f(x)}M$ so that $f(x) = \varphi_f(x) + \sigma_f(x)$ with . In the following we view $\sigma_f$ as a section of the pull--back bundle $\varphi_f^*\nu M$. Observe that $\sigma_f = 0$ if and only if $f(M) \subset M$.

\medskip

The operator $P_t$ considered as acting on $\R^n$--valued functions on $M$ is
not elliptic as the computation of its symbol shows below. In order to complement it to an elliptic operator we proceed as follows: For $\varphi \in \diff(M)_+$
consider the connection Laplacian $\Delta^{\varphi^*\nu M} =
(\nabla^{\varphi^*\nu M})^*\nabla^{\varphi^*\nu M}$ on the pull--back bundle
$\varphi^*\nu M$. Here $\nabla^{\varphi^*\nu M}$ denotes the pull--back
connection on $\varphi^* \nu M$. Its formal adjoint is taken with respect
to the metric $g_{\Omega_0}$ on $M$ and the pull--back metric $\varphi^* h$ on
$\varphi^* \nu M$. Then, reasoning similarly as above, $f \mapsto \Delta^{\varphi_f^*\nu M}\sigma_f:\mathcal
U \subset C^{\infty}(M,\R^n)\to C^{\infty}(M, \varphi_f^*\nu M) \subset C^{\infty}(M,\R^n)$ is a quasilinear, second order differential operator, and so is
\ben
\widetilde P_t:\mc{U}\subset C^{\infty}(M,\R^n)\to C^{\infty}(M,\R^n),\quad f\mapsto P_t(f)-\Delta^{\varphi_f^*\nu M}\sigma_f.
\ee
We wish to establish short--time existence and uniqueness of the associated flow equation
\be\label{auxflow}
\frac{\partial}{\partial t}f(t)=\widetilde P_t\big(f(t)\big),\quad f(0)=f_0.
\ee
To compute the linearisation $D_{f_0}\widetilde P_0(Y)$ we take a curve $f_s\subset\mc{U}$ through $f_0$ with
\ben
Y(x)=\frac{d}{ds}f_s(x)|_{s=0}\in T_x\R^n\cong\R^n.
\ee 
We write $Y^\parallel(x)$ and $Y^\perp(x)$ for the projections of $Y(x)$ to $T_xM$ and $\nu_xM$. By design of the extension of $\Omega_0$ to $\R^n$ (cf.\ our convention above), $\varphi_s:=\pi\circ f_s\in\diff(M)_+$ satisfies $f^*_s\Omega_0=\varphi^*_s\Omega_0$. Furthermore, one clearly has
\be\label{tangent}
\frac{d}{ds}\varphi_s(x)|_{s=0}=Y^\parallel(x) \quad \text{and} \quad \frac{d}{ds}\sigma_s(x)|_{s=0} = Y^\perp(x)  
\ee
for $\sigma_s(x):=f_s(x)-\varphi_s(x) \in \nu_{\varphi_s(x)}M$.

\medskip

First we compute the linearisation of $P_0$. Using the naturality of the vector field $X_{\Omega_0}$, i.e.~$\varphi_*X_{\varphi^*\Omega_0}\big(\Omega(t)\big)=X_{\Omega_0}\big(\varphi^{-1*}\Omega(t)\big)$ for $\varphi \in \diff(M)_+$, we obtain
\begin{eqnarray*}
\frac{d}{ds}P_0(f_s)|_{s=0} &=& -\frac{d}{ds}df_s\big(X_{\varphi^*_s\Omega_0}(\Omega_0)\big)|_{s=0}\\
 &=& -\frac{d}{ds} df_s \big( \varphi_{s*}^{-1} X_{\Omega_0} ( \varphi_s^{-1*}\Omega_0)\big)|_{s=0}.
\end{eqnarray*}
For $Y$ tangent to $M$, i.e.~$Y(x)=Y^\parallel(x)$ for all $x\in M$, we may in view of \eqref{tangent} assume that $f_s=\varphi_s \in \diff(M)_+$ for all $s$ and we get
\begin{eqnarray*}
\frac{d}{ds}P_0(f_s)|_{s=0} &=& -\frac{d}{ds} X_{\Omega_0}(\varphi^{-1*}_{s}\Omega_0)|_{s=0}\\
&=& X_{\Omega_0}(\mc{L}_Y\Omega_0)\\ 
&=& -(\delta_{\Omega_0}d(Y\llcorner\Omega_0))\llcorner\Omega_0+\mbox{ terms of lower order in }Y. 
\end{eqnarray*}
For $Y$ perpendicular to $M$, i.e.~$Y(x)=Y^\perp(x)$ for all $x \in M$, we may assume that $\varphi_s=\id$ for all $s$, again using \eqref{tangent}. Hence
\ben
\frac{d}{ds}P_0(f_s)|_{s=0} =  -\frac{d}{ds} df_s X_{\Omega_0}(\Omega_0) = dY(X_{\Omega_0}(\Omega_0)),
\ee
which is of lower order in $Y$ and hence does not contribute to the symbol. For general $Y=Y^\parallel + Y^\perp$ we therefore find
\ben
D_{f_0} P_0 (Y) =  -(\delta_{\Omega_0}d(Y^\parallel\llcorner\Omega_0))\llcorner\Omega_0+\mbox{ terms of lower order in }Y.
\ee
For the linearisation of $f \mapsto \Delta^{\varphi_f^*\nu M}\sigma_f$ we again assume $Y$ to be tangent to $M$ first. For a curve $f_s=\varphi_s \in \diff(M)_+$ as above we have in particular that $\sigma_s=0$ for all $s$. Hence $\frac{d}{ds}\Delta^{\varphi_s^*\nu M}\sigma_s|_{s=0}=0$. For $Y$ perpendicular to $M$ we have $\sigma_s(x) = f_s(x)-x$ for a curve $f_s$ with $\varphi_s \equiv \id$ as above. Hence $\frac{d}{ds} \Delta^{\varphi_s^*\nu M}\sigma_s|_{s=0}=\Delta^{\nu M} Y$. For general $Y=Y^\parallel+Y^\perp$ we get
\ben
D_{f_0}(f \mapsto \Delta^{\varphi_f^*\nu M}\sigma_f)(Y)=\Delta^{\nu M}Y^\perp.
\ee
It follows that the symbol of the linearised operator $D_{f_0}\widetilde P_0$ is
\ben
\sigma(D_{f_0}\widetilde P_0)(x,\xi)Y=-\big(\xi\llcorner\big(\xi\wedge(Y^\parallel\llcorner\Omega_0)\big)\big)\llcorner\Omega_0-|\xi|^2_{\Omega_0}Y^\perp.
\ee
To check strong parabolicity we assume $|\xi|_{\Omega_0}=1$ and write
$Y^\parallel=a\xi+Y_0$, $a\in\R$, $Y_0\in\xi^\perp$ and $\Omega_0=\omega\wedge\xi+\psi_+$. Then
\begin{eqnarray*}
-(\pi^*g_{\Omega_0}+h)\big(\sigma(D_{f_0}\widetilde P_0)(x,\xi)Y,Y\big) & = & g_{\Omega_0}\big(\big(\xi\llcorner\big(\xi\wedge\big(Y^\parallel\llcorner\Omega_0\big)\big)\big)\llcorner\Omega_0,Y^\parallel\big)+|Y^\perp|^2_h\\
& = & g_{\Omega_0}\big(\Omega_0,\big(\xi\llcorner\big(\xi\wedge(Y^\parallel\llcorner\Omega_0)\big)\big)\wedge Y^\parallel\big)+|Y^\perp|^2_h\\
& = & g_{\Omega_0}\big(Y^\parallel\llcorner\Omega_0,\xi\llcorner\big(\xi\wedge(Y^\parallel\llcorner\Omega_0)\big)\big)+|Y^\perp|^2_h\\
& = & g_{\Omega_0}\big(a\omega+Y_0\llcorner\psi_++(Y_0\llcorner\omega)\wedge\xi,a\omega+Y_0\llcorner\psi_+\big)\\
& & +|Y^\perp|^2_h\\
& = & 3|a|^2_{\Omega_0}+|Y_0\llcorner\psi_+|^2_{\Omega_0}+|Y^\perp|^2_h\\
& \stackrel{\eqref{su3isom}}{=} & 3|a|^2+2|Y_0|^2+|Y^\perp|^2_h\\
& \geq & (\pi^*g_{\Omega_0}+h)(Y,Y).
\end{eqnarray*}
Theorem~\ref{mainanaresult} applies once again to yield short--time existence
and uniqueness of~\eqref{auxflow}. Last we show that a solution $f(t)$
to~\eqref{auxflow} satisfies $f(t)(M)\subset M$ for all $t$. As in this case
$\sigma_{f(t)}$ is just the zero section of $\varphi(t)^*\nu M$, we obtain the
desired solution to~\eqref{unicflow}. We proceed as in the harmonic map heat
flow case, see for instance Part IV in \cite{ha75}: Consider the bundle endomorphism $r:\nu M\to\nu M$ which is
multiplication by $-1$ in each fibre. We claim that $\widetilde P_t\circ
r=dr\circ\widetilde P_t$. This is clear for $P_t$ since $(r\circ
f)^*\Omega_0=f^*\Omega_0$. Furthermore, $\varphi_{r \circ f} = \varphi_f$ and $\sigma_{r \circ f} = r \circ \sigma_f$, such that by linearity of the connection Laplacian one has 
$$
\Delta^{\varphi_{r \circ f}^*\nu M} \sigma_{r \circ f} = r \circ \Delta^{\varphi_f^*\nu M} \sigma_f = dr \circ \Delta^{\varphi_f^*\nu M} \sigma_f,
$$
where the last equality follows by viewing a section of $\nu M$ as a vertical
vector field on the total space of that bundle. This proves our claim and implies that for a solution $f(t)$ to~\eqref{auxflow} the composition $r \circ
f(t)$ is again a solution with initial condition $f_0$. Now if $f(t)(M)$ were
not contained in $M$ for some $t$, then $r\circ f$ would yield a second, different solution with same initial condition, contradicting uniqueness.
\end{prf}

A solution $f(t)\in\mc{U}$ to~\eqref{finalflow} yields a solution
$\varphi_t=f_0^{-1}\circ f(t)\in\diff(M)_+$ to~\eqref{backwardflow} for a
given Dirichlet energy flow solution $\Omega(t)$. From there, uniqueness easily follows:

\begin{cor}\label{DF-unique}
Suppose that $\Omega(t)$ and $\Omega'(t)$ are two solutions to~\eqref{floweq} for $t\in[0,\epsilon]$, $\epsilon>0$. If $\Omega(0)=\Omega_0=\Omega'(0)$, then $\Omega(t)=\Omega'(t)$ for all $t\in[0,\epsilon]$. 
\end{cor}
\begin{prf}
Solving for~\eqref{backwardflow} with $\Omega(t)$ and $\Omega'(t)$ gives two flows $\varphi_t$ and $\varphi'_t$ which without loss of generality we assume to be defined on $[0,\epsilon]$. By design $\widetilde\Omega(t)=\varphi^*_t\Omega(t)$ and $\widetilde\Omega'(t)=\varphi'^*_t\Omega'(t)$ define a solution to~\eqref{ddflow} at $\Omega_0$. Uniqueness of the Dirichlet--DeTurck flow implies $\widetilde\Omega(t)=\widetilde\Omega'(t)$. Hence $\varphi_t$ and $\varphi'_t$ are solutions of the ordinary differential equation
\be\label{finalode}
\frac{\partial}{\partial t}\psi_t=-X_{\Omega_0}\big(\widetilde\Omega(t)\big)\circ\psi_t.
\ee
By uniqueness of the solution to~\eqref{finalode}, we conclude $\varphi_t=\varphi'_t$, whence $\Omega(t)=\Omega'(t)$.
\end{prf}
%
%
%
%
%
\section{The second variation of $\mc{D}$}
\label{linstab}
In this section we compute the second variation of $\mc{D}$ at some fixed $\bar\Omega\in\mc{X}=Q^{-1}(0)$ (cf.\ Corollary~\ref{critker}). Further, we show that $\mc{X}$ is a Fr\'echet manifold whose tangent space at $\bar\Omega$ is precisely $\ker D^2_{\bar\Omega}\mc{D}$. All metric quantities, projections on $\Gt$--invariant modules etc.~will be taken \wrt this $\bar\Omega$.

\medskip

For a given vector bundle $E\to M$ we denote by $W^{k,2}(E)$ the space of sections whose local components have square integrable derivatives up to order $k$. The associated Sobolev norm will be written $\|\cdot\|_{W^{k,2}}$. Further, we simply write $L^2(E)$ for $W^{0,2}(E)$. More generally, we consider the Hilbert manifolds $W^{k,2}(\xi)$ for a fibre bundle $\xi\to M$ in order to deal with non--linear differential operators. The integer $k$ will be chosen appropriately when required, but at any rate big enough so that all sections involved are at least of class $C^0$ and the corresponding function spaces $W^{k,2}(M,\R)$ are Banach algebras under pointwise multiplication. According to the Sobolev embedding theorem we therefore need $k>11/2$ for $Q$ to extend to a smooth map
\ben
Q_k:W^{k,2}\big(\Lambda^3_+T^*M\big)\to W^{k-2,2}\big(\Lambda^3T^*M\big).
\ee

\begin{lem}\label{xspace}
The space $\mc{X}^k:=Q_k^{-1}(0)$ is a smooth Banach manifold whose tangent space at $\bar\Omega\in\mc{X}$ is given by
\ben
T_{\bar \Omega}\mc{X}^k=\{\dot\Omega\in W^{k,2}(\Lambda^3T^*M)\,|\,d\dot\Omega=0,\,d\dot\Theta_{\bar\Omega}=0\}.
\ee
\end{lem}
\begin{prf}
We have $\mc{X}^k=\{\Omega\in W^{k,2}(\Lambda^3_+T^*M)\,|\,d\Omega=0,\,d\Theta(\Omega)=0\}$, so it remains to show that the extension $N_k$ of
\ben
N:\Omega\in\Omega^3_+(M)\mapsto\big(d\Omega,d\Theta(\Omega)\big)\in\Omega^4(M)\times\Omega^5(M)
\ee
to $W^{k,2}(\Lambda^3_+T^*M)$ has $0$ as a regular value. By the first example of Section~\ref{def}, $(\dot N_k)_{\bar\Omega}=\big(d\dot\Omega,d\!\star\! p(\dot\Omega)\big)$. Since the range of $(d_p)_k:W^{k,2}(\Lambda^pT^*M)\to W^{k-1,2}(\Lambda^{p+1}T^*M)$ is closed, $(B^{p+1})^{k-1}=\im (d_p)_k$ is a Banach space. Hence for a given $(d\sigma,d\tau)\in(B^4)^{k-1}\times (B^5)^{k-1}$, we need $\dot\Omega\in\Omega^3(M)$ \st $d\dot\Omega=d\sigma$ and $d\!\star\! p(\dot\Omega)=d\tau$. The Hodge decomposition theorem of $[\dot\Omega]_q$ gives $\mc{H}(\dot\Omega_q)\oplus d\dot\alpha_q\oplus\delta\dot\beta_q$, where $\dot\alpha_q\in W^{k+1,2}(\Lambda^2T^*M)$, $\dot\beta_q\in W^{k+1,2}(\Lambda^4T^*M)$ for $q=1,\,7,\,27$, and $\mc{H}$ denotes projection on the space of harmonic forms. Similar decompositions hold for $\sigma$ and $\tau$. Taking $\dot\beta_q=[\beta_\sigma]_q$ yields
\ben
d\dot\Omega=\bigoplus_{q\in\{1,\,7,\,27\}}d\delta\dot\beta_q=d\delta\beta_\sigma=d\sigma.
\ee
On the other hand,
\ben
d\!\star\! p(\dot\Omega)=d\!\star\!\big(\tfrac{4}{3}d\dot\alpha_1+d\dot\alpha_7-d\dot\alpha_{27}\big)=d\delta\beta_\tau=d\tau,
\ee
provided we put $\dot\alpha_1=3\star[\beta_\tau]_1/4$, $\dot\alpha_7=\star[\beta_\tau]_7$ and $\dot\alpha_{27}=\star[\beta_\tau]_{27}$. Consequently, $(\dot N_k)_{\bar\Omega}$ is surjective, whence the result by the Banach space implicit function theorem.
\end{prf}

Consider the closed linear subspace $V_{\bar\Omega}^k=\{\dot\Omega\in W^{k,2}(\Lambda^3T^*M)\,|\,[\delta\dot\Omega]_7=0\}\subset\lambda^{-1}(0)$ (cf.~\eqref{lambdaomega}) and let
\ben
S_{\bar\Omega}:=V_{\bar\Omega}^k\cap\mc{X}^k.
\ee
Note that $S_{\bar\Omega}\subset\widetilde Q_{\bar\Omega}^{-1}(0)$, whence $S_{\bar\Omega}\subset\Omega^3(M)$ for $\widetilde Q_{\bar\Omega}$ is a quasilinear, elliptic operator by Proposition~\ref{qtildeoperator}.

\begin{prp}\label{slice}
Near $\bar\Omega$, the space $S_{\bar\Omega}$ is a smooth submanifold of $\mc{X}^k$. Its tangent space at $\bar\Omega$ is naturally isomorphic with the space of $\bar \Omega$--harmonic $3$--forms. In particular, $\dim S_{\bar\Omega}=b_3$. 
\end{prp}

\begin{rem}
In Section~\ref{stabsec} we actually show that near $\bar\Omega$, $\widetilde Q_{\bar\Omega}^{-1}(0)$ coincides with $S_{\bar\Omega}$ (cf.\ Corollary~\ref{zero_set}).
\end{rem}

Before we can prove Proposition~\ref{slice}, we need a technical result first. Consider the Hilbert manifold $\diff(M)_0^{k+1}$ obtained from completion of the identity component of $\diff(M)$ with respect to $\|\cdot\|_{W^{k+1,2}}$. Then $\diff(M)_0^{k+1}$ is a smooth Banach manifold and a topological group which acts continuously on $\mc{X}^k$ via pull--back. We denote by $\mr{I}_{\bar\Omega}$ the subgroup of $\diff(M)^{k+1}_0$ fixing $\bar \Omega$, i.e.\ $\varphi^*\bar \Omega=\bar \Omega$ for all $\varphi\in\mr{I}_{\bar\Omega}$; in particular $g=\varphi^*g$ which by~\cite{ta06} implies $\varphi\in\diff(M)_0$ since any $\varphi\in\diff(M)^{k+1}_0$ is of class $C^1$. Hence $\mr{I}_{\bar\Omega}$ is contained in $\diff(M)_0$. Let $\mc{O}^k_{\bar\Omega}$ denote the $\diff(M)_0^{k+1}$--orbit through $\bar\Omega$ in $\mc{X}^k$ and put
\ben
G_{\bar\Omega}^{k+1}:=\diff(M)^{k+1}_0/\mr{I}_{\bar\Omega}.
\ee
As in Section 5 of~\cite{eb68} one can prove that $G_{\bar\Omega}^{k+1}$ endowed with the quotient topology is a smooth manifold.

\begin{lem}\label{orbit}
If $k>11/2$, then $[\varphi]\in G_{\bar\Omega}^{k+1}\mapsto\varphi^*\bar\Omega\in\mc{X}^k$ is an injective immersion with closed image. In particular, $\mc{O}^k_{\bar\Omega}$ is a closed, smooth submanifold of $\mc{X}^k$ with tangent space
\ben
T_{\bar \Omega}\mc{O}^k_{\bar\Omega}=\{\mc{L}_X\bar \Omega\,|\,X\in W^{k+1,2}(TM)\}.
\ee
\end{lem}
\begin{prf}
We can argue as in Ebin's proof of the corresponding result for the moduli space of Riemannian metrics, cf.\ Section 6 in~\cite{eb68}. The only remaining issue to check is the injectivity of the symbol of $\lambda^*:C^{\infty}(TM)\to\Omega^3(M)$ defined in~\eqref{lambdastern} (this ensures that the extension to a map $W^{k+1,2}(TM)\to W^{k,2}(\Lambda^3T^*M)$ has closed range). Indeed,
\be\label{injsym}
\sigma(\lambda^*)(x,\xi)(v)=i\xi\wedge(v\llcorner\bar \Omega),
\ee
and this vanishes \iff $v\llcorner\bar \Omega$ is of the form $\eta\wedge\xi$ for some $\eta\in\Omega^1(M)$. But $v \llcorner \bar\Omega \in \Lambda^2_{7}$, so that in this case $0=[\eta\wedge\xi]_{14}=(2\eta\wedge\xi-\star(\eta\wedge\xi\wedge\Omega))/3$. Since on the right hand side the first
term contains $\xi$ while the second does not, this can only hold if $\eta\wedge\xi=0$, i.e.~if $\eta=0$ or equivalently, $v=0$.
\end{prf}

\begin{prf} (of Proposition~\ref{slice}) Since $\mc{O}^k_{\bar\Omega}\subset\mc{X}^k$, the tangent space $T_{\bar\Omega}\mc{O}^k_{\bar\Omega}$ is contained in $T_{\bar\Omega}\mc{X}^k$. By extending Proposition~\ref{divsplit} to Sobolev spaces, we deduce $T_{\bar\Omega}\mc{X}^k+T_{\bar\Omega}V_{\bar\Omega}^k=W^{k,2}(\Lambda^3T^*M)$. Hence, the intersection is transversal near $\bar\Omega$ so that $S_{\bar\Omega}$ is a smooth submanifold of $W^{k,2}(\Lambda^3T^*M)$ in a neighbourhood of $\bar\Omega$. The tangent space at $\bar\Omega$ is 
\ben
T_{\bar\Omega}S_{\bar\Omega}=T_{\bar\Omega} V_{\bar\Omega}^k\cap T_{\bar\Omega}\mc{X}^k=\{\dot\Omega\in\Omega^3(M)\,|\,d\dot\Omega=0,\,d\dot\Theta_{\bar\Omega}=0,\,[\delta\dot\Omega]_7=0\}.
\ee
Next, the map $\dot\Omega\in T_{\bar\Omega}S_{\bar\Omega}\mapsto[\dot\Omega]\in H^3(M,\R)$ is an isomorphism. For injectivity, assume $\dot\Omega=d\eta$ so that $[\delta d\eta]_7=0$. By Lemma 10.3.2 in~\cite{jo00}, this implies $\delta[d\eta]_1=0$ and $\delta[d\eta]_7=0$, whence $\delta[d\eta]_{27}=0$ for $d\!\star\! p(d\eta)=0$. Consequently, $\dot\Omega=d\eta$ is $\Delta$--harmonic which is impossible unless $\dot\Omega=0$. For surjectivity, recall that the projections on irreducible $\Gt$--components $\Lambda^p_ qT^*M$ commute with the Hodge Laplacian $\Delta$ since $\bar\Omega$ is torsion--free (cf.\ for instance Theorem 3.5.3 in~\cite{jo00}). Hence, a $p$--form is $\Delta$--harmonic \iff its irreducible components in $\Omega^p_q(M)$ are $\Delta$--harmonic. So, given a cohomology class $c=[\dot\Omega]\in H^3(M,\R)$ with unique $\Delta$--harmonic representative $\dot\Omega$, we have $d\!\star\! p(\dot\Omega)=0$ and thus $\dot\Omega\in T_{\bar\Omega}S_{\bar\Omega}$.
\end{prf}

In particular, we deduce in conjunction with Proposition~\ref{divsplit} that
\be\label{transversality}
T_{\bar\Omega}\mc{X}^k=T_{\bar\Omega}\mc{O}^k_{\bar\Omega}\oplus T_{\bar\Omega}S_{\bar\Omega}.
\ee
Next consider the map
\ben
\Phi^k:([\varphi],\Omega)\in G^{k+1}_{\bar\Omega}\times\mc{S}_{\bar\Omega}\mapsto[\varphi]^*\Omega\in\mc{X}^k
\ee
on a suitable neighbourhood $\mc{S}_{\bar\Omega}$ of $\bar\Omega$ in $S_{\bar\Omega}$. Following~\cite{jo00} $\Phi^k$ is well--defined: Indeed, near $\bar\Omega$ we can linearise the action of $\mr{I}_{\bar\Omega}$ on $S_{\bar\Omega}$ via the exponential map $\exp_{\bar\Omega}:T_{\bar\Omega}S_{\bar\Omega}\cong H^3(M,\R)\to S_{\bar\Omega}$ induced by $\langle\cdot\,,\cdot\rangle_{L^2}$. The linearised action must be trivial as any element in $\mr{I}_{\bar\Omega}$ is homotopic to the identity. Hence $\mr{I}_{\bar\Omega}$ acts trivially on $S_{\bar\Omega}$ close to $\bar\Omega$. The Banach space inverse function theorem and~\eqref{transversality} imply that $\Phi^k$ is a diffeomorphism onto its image near $([\Id_M],\bar\Omega)$. Hence, shrinking $\mc{S}_{\bar\Omega}$ possibly further, we have shown that $\mc{S}_{\bar\Omega}$ defines a {\em slice} for the $\diff(M)^{k+1}_0$--action on $\mc{X}^k$ near $\bar\Omega$. The same statement holds for the $C^\infty$--topology instead of the Sobolev topologies, whence the

\begin{cor}\label{xmanif}
(i) The space $\mc{X}$ is a Fr\'echet manifold whose tangent space at $\bar\Omega$ is given by
\ben
T_{\bar\Omega}\mc{X}=\{\dot\Omega\in\Omega^3(M)\,|\,d\dot\Omega=0,\,d\dot\Theta_{\bar\Omega}=0\}.
\ee
(ii) {\bf (Joyce)} The space of torsion--free $\Gt$--structures modulo diffeomorphisms isotopic to the identity is a smooth manifold of dimension $b_3=\dim H^3(M,\R)$.
\end{cor}

We can now prove the central result of this section.

\begin{prp}\label{secvar}
Let $\bar\Omega\in\mc{X}$.

(i) We have
\ben
D^2_{\bar\Omega}\mc{D}(\dot\Omega,\dot\Omega)=\int_M\big(|d\dot\Omega|^2+|d\dot\Theta_{\bar\Omega}|^2\big)\vol_{\bar\Omega}.
\ee
In particular, the second variation of $\mc{D}$ at $\bar\Omega$ is a positive semi--definite bilinear form with
\ben
\ker D^2_{\bar\Omega}\mc{D}=T_{\bar\Omega}\mc{X}.
\ee

(ii) The linearisation $L_{\bar\Omega}:=D_{\bar\Omega}\widetilde Q_{\bar\Omega}$ of $\widetilde Q_{\bar\Omega}$ is a symmetric, non--positive and elliptic operator given by
\ben
L_{\bar\Omega}\dot\Omega=-\delta d\dot\Omega-pd\delta p\dot\Omega-3d[\delta\dot\Omega]_7.
\ee
More precisely, writing $\dot\Omega=\dot f\bar\Omega\oplus\star(\dot\alpha\wedge\bar\Omega)\oplus\dot\gamma$, we have
\be\label{pseudowb}
-L_{\bar\Omega}\dot\Omega=\Delta\dot\Omega+\tfrac{34}{21}d^7_1d^1_7\dot
f\cdot\bar\Omega+\star\big(d^7_7d^7_7\dot\alpha\wedge\bar\Omega\big)+d^7_{27}d^{27}_7\dot\gamma-\tfrac{2}{21}d^7_1d^{27}_7\dot\gamma\cdot\bar\Omega-\tfrac{2}{3}d^7_{27}d^1_7\dot f
\ee
for the $\Gt$--differential operators $d^p_q$ introduced in Section~\ref{reptheory}. In particular
\ben
\langle L_{\bar\Omega}\dot\Omega,\dot\Omega\rangle_{L^2_{\bar\Omega}}=-\|d\dot\Omega\|_{L^2_{\bar\Omega}}^2-\|\delta p\dot\Omega\|_{L^2_{\bar\Omega}}^2-3\|[\delta\dot\Omega]_7\|_{L^2_{\bar\Omega}}^2\mbox{ and }\ker L_{\bar\Omega}=T_{\bar\Omega}\mc{S}_{\bar\Omega}.
\ee 
\end{prp}
\begin{prf}
(i) We compute the second variation of $\mathcal D$ at a critical point $\bar \Omega$ by differentiating equation \eqref{gradD} once more. Thus we obtain
\ben
D^2_{\bar\Omega}\mc{D}(\dot\Omega,\dot\Omega)=\int_Md\dot\Omega\wedge\star d\dot\Omega+d\dot\Theta_{\bar\Omega}\wedge\star d\dot\Theta_{\bar\Omega}=\int_M\big(|d\dot\Omega|^2+|d\dot\Theta_{\bar\Omega}|^2\big)\vol \geq 0,
\ee
for the remaining terms involve either $d\bar\Omega$ or $d\Theta(\bar\Omega)$. But these terms vanish since $\bar\Omega$ is a critical point and hence torsion--free. Furthermore, the kernel $D^2_{\bar\Omega}\mc{D}$ is precisely $T_{\bar\Omega}\mc{X}$ by Corollary~\ref{xmanif}.

\medskip

(ii) Since $D^2_{\bar\Omega}\mc{D}(\dot\Omega,\dot\Omega)=\langle\mr{Hess}_{\bar\Omega}\mc{D}\dot\Omega,\dot\Omega\rangle_{L^2}$ and $D_{\bar\Omega}Q=-\mr{Hess}_{\bar\Omega}\mc{D}$, it follows from (i) that $D_{\bar\Omega} Q=-\delta d\dot\Omega-p d\delta p\dot\Omega$. On the other hand, using again that $\bar \Omega$ is torsion--free, we get
\ben
D_{\bar\Omega}\Lambda_{\bar\Omega}(\dot\Omega)=d\big(X_{\bar\Omega}(\dot\Omega)\llcorner\bar\Omega\big).
\ee
Bearing~\eqref{7contrOm} in mind, we obtain $\dot\Lambda_{\bar\Omega}=-3d[\delta\dot\Omega]_7$, whence
\ben
L_{\bar\Omega}(\dot\Omega)=-\delta d\dot\Omega-pd\delta p\dot\Omega-3d[\delta\dot\Omega]_7.
\ee
This operator is clearly symmetric and non--positive. Ellipticity was proven in Lemma~\ref{qtildeoperator}. To compute~\eqref{pseudowb} we start from~\eqref{deltapq} and use Tables~\ref{extderfor} and~\ref{secordid}. Then
\begin{eqnarray*}
\delta d\dot\Omega & = & \tfrac{1}{7}(4d^7_1d_7^1\dot f+d^7_1d^{27}_7\dot\gamma)\cdot\bar\Omega\oplus\star\big((d^1_7d^7_1\dot\alpha+\tfrac{1}{2}d^7_7d^7_7\dot\alpha+\tfrac{1}{4}d^7_7d^{27}_7\dot\gamma)\wedge\bar\Omega\big)\\
& & \oplus\, d^7_{27}d^1_7\dot f+d^7_{27}d^7_7\dot\alpha+\frac{1}{4}d^7_{27}d^{27}_7\dot\gamma+d^{27}_{27}d^{27}_{27}\dot\gamma\\
d\delta\dot\Omega & = & \tfrac{1}{7}(3d^7_1d_7^1\dot f-d^7_1d^{27}_7\dot\gamma)\cdot\bar\Omega\oplus\star\big(\tfrac{1}{2}\big(d^7_7d^7_7\dot\alpha-\tfrac{1}{2}d^7_7d^{27}_7\dot\gamma)\wedge\bar\Omega\big)\\
& & \oplus-d^7_{27}d^1_7\dot f-d^7_{27}d^7_7\dot\alpha+\tfrac{1}{3}d^7_{27}d^{27}_7\dot\gamma+d^{14}_{27}d^{27}_{14}\dot\gamma
\end{eqnarray*}
from which~\eqref{pseudowb} follows by applying~Table~\ref{laplace}.
\end{prf}

Let $\Db$ denote the associated Hodge--Dirac operator \wrt the metric induced by $\bar\Omega\in\mc{X}$, i.e.
\ben
\Db\dot\Omega=d\dot\Omega+\delta\dot\Omega.
\ee
In particular, $\Db$ is symmetric and $\Db^2=\Delta$. In view of longtime existence to be established in the next section we note the following corollary.

\begin{cor}[G\r arding inequality]\label{coercivity}
For all $\dot\Omega\in\Omega^3(M)$, 
\end{cor}
\ben
\langle-L_{\bar\Omega}\dot\Omega,\dot\Omega\rangle_{L^2}\geq \|\Db \dot\Omega\|^2_{L^2}.
\ee
{\it In particular, we have
\ben
\langle-L_{\bar\Omega}\dot\Omega,\dot\Omega\rangle_{L^2}\geq C\|\dot\Omega\|^2_{W^{1,2}}-\|\dot\Omega\|^2_{L^2}
\ee
for some constant $C$ independent of $\dot\Omega$.}
\begin{prf}
Writing $\dot\Omega$ as in Proposition~\ref{secvar}, the first inequality follows from
\begin{eqnarray}\label{quad_form_L}
\langle-L_{\bar\Omega}\dot\Omega,\dot\Omega\rangle_{L^2} & = & \|\Db\dot\Omega\|^2_{L^2}+\tfrac{34}{3}\|d^1_7\dot f\|^2_{L^2}+4\|d^7_7\dot\alpha\|^2_{L^2}+\|d^{27}_7\dot\gamma\|^2_{L^2}-\tfrac{4}{3}\langle d^{27}_7\dot\gamma,d^1_7\dot f\rangle_{L^2}\\
& \geq &  \|\Db \dot\Omega\|^2_{L^2}+\tfrac{34}{3}\|d^1_7\dot f\|^2_{L^2}+\|d^{27}_7\dot\gamma\|^2_{L^2}-\tfrac{2}{3}(\|d^1_7\dot f\|^2_{L^2}+\|d^{27}_7\dot\gamma\|^2_{L^2})\nonumber\\
& \geq & \|\Db \dot\Omega\|^2_{L^2}.\nonumber
\end{eqnarray}
The second inequality is just the elliptic estimate $\|\dot\Omega\|^2_{W^{1,2}}\leq C^{-1}(\|\dot\Omega\|^2_{L^2}+\|\Db\dot\Omega\|^2_{L^2})$ for some constant $C^{-1}$.
\end{prf}
%
%
%
%
%
\section{Stability}\label{stabsec}
We continue to fix a torsion--free $\Gt$--form $\bar\Omega\in\mc{X}$. From now on, we let $W^{k,2}$ and $W^{k,2}_+$ be shorthand for the Sobolev spaces $W^{k,2}(\Lambda^3T^*M)$ and $W^{k,2}(\Lambda^3_+T^*M)$ \wrt the metric $g_{\bar\Omega}$. The induced norm will be denoted by $\|\cdot\|_{W^{k,2}}$. Again, $k$ is an integer strictly greater than $11/2$ which for simplicity we assume to be odd (this avoids using fractional Sobolev spaces below). In particular, $W^{k,2}$ embeds continuously into $C^2$. The goal of this section is to prove the subsequent stability theorem.

\begin{thm}[Stability]\label{stability} Let $\bar\Omega\in\Omega^3_+(M)$ be a torsion--free $\Gt$--form. For all $\epsilon>0$ there exists some $\delta=\delta(\epsilon)>0$ such that for any $\widetilde\Omega_0$ with $\|\widetilde\Omega_0-\bar\Omega\|_{W^{k,2}}<\delta$, the Dirichlet--DeTurck flow $\widetilde\Omega_t$ at $\bar\Omega$ with initial condition $\widetilde\Omega_0$  
\begin{enumerate}
\item {\bf (longtime existence)} exists for all $t\in[0,\infty)$, 
\item {\bf (a priori estimate)} satisfies the estimate $\|\widetilde\Omega_t-\bar\Omega\|_{W^{k,2}}<\epsilon$ for all $t\in[0,\infty)$, and
\item {\bf (convergence)} converges \wrt the $W^{k,2}$--norm to a torsion--free $\Gt$--form $\widetilde{\Omega}_\infty$ as $t\rightarrow\infty$.
\end{enumerate}
\end{thm}

Since the Dirichlet energy flow exists as long as the Dirichlet--DeTurck flow exists we immediately obtain:

\begin{cor}\label{convcor}
Let $\bar\Omega\in\Omega^3_+(M)$ be a torsion--free $\Gt$--form. For initial
conditions sufficiently  $C^\infty$--close to $\bar\Omega$ the Dirichlet
energy flow $\Omega_t$ exists for all times and converges modulo diffeomorphisms to a
torsion--free $\Gt$-form $\Omega_\infty$, i.e.~there exists a family of diffeomorphisms $\varphi_t \in
\diff(M)_+$ such that $\varphi_t^*\Omega_t$ converges to $\Omega_\infty$ with
respect to the $C^\infty$--topology.
\end{cor}

The proof of Theorem \ref{stability} will be subdivided into a sequence of intermediate steps.

\bigskip

First, we tackle existence of the Dirichlet--DeTurck flow together with the a priori estimate on arbitrary, but finite time intervals for initial conditions sufficiently close to $\bar
\Omega$. Here we use the Banach space inverse function theorem, following the
approach of Huisken and Polden~\cite{hp99} for geometric evolution equations
for hypersurfaces. Let $0<T<\infty$. If $\pi:M\times[0,T]\to M$ denotes projection onto the first factor, let $C^\infty(M \times [0,T],\pi^*\Lambda^3)$ denote the space of smooth sections of $\Lambda^3$ pulled back to $M\times[0,T]$. For any non--negative integer $s$ we define the Hilbert space $V^s[0,T]$ as the completion of $C^\infty(M \times [0,T],\pi^*\Lambda^3)$ with respect to the inner product given by
\ben
\langle\dot\Omega_1,\dot\Omega_2\rangle_{V^s[0,T]}=\sum_{j\leq s}\int_0^Te^{-2t}\langle\partial^j_t\dot\Omega_1,\partial^j_t\dot\Omega_2\rangle_{W^{2(s-j),2}}dt.
\ee 
In particular, $\|\dot\Omega\|^2_{V^0[0,T]}=\int_0^Te^{-2t}\|\dot\Omega\|^2_{L^2}dt$. These spaces are also known as {\em anisotropic} or {\em parabolic} Sobolev spaces, where a time derivative has the weight of two space derivatives. The density $e^{-2t}$ is introduced for technical reasons, see below. Similarly, we can consider the Hilbert manifold $V^s_+[0,T]$ consisting of sections of $\pi^*\Lambda^3_+$ of class $V^s$. Define the map
\be\label{banachflowmap}
F:V^s_+[0,T]\to W^{2s-1,2}_+\times V^{s-1}[0,T],\quad F(\widetilde\Omega)=\bigl(\widetilde\Omega_0,\partial_ t\widetilde\Omega_t-\widetilde Q _{\bar\Omega}(\widetilde\Omega_t)\bigr).
\ee
As usual, restricting to the boundary is tantamount to invoking a trace theorem, which in this context is stating that the trace map $\widetilde\Omega\mapsto\widetilde\Omega_0$ is continuous from $V^s[0,T]$ to $W^{2s-1,2}$, see \cite{hp99}. We wish to show that $F$ is a local diffeomorphism near $\bar\Omega$. Consider the linearisation at $\bar\Omega$ of the map~\eqref{banachflowmap}, namely
\be\label{banachflowlin}
D_{\bar\Omega}F:V^s[0,T]\to W^{2s-1,2}\times V^{s-1}[0,T],\quad
D_{\bar\Omega}F(\dot\Omega)=(\dot\Omega_0,P\dot\Omega),
\ee
where $P\dot\Omega:=\partial_t\dot\Omega_t-L_{\bar\Omega}\dot\Omega_t$. Let $\mb H$ be the completion of $C^\infty(M \times [0,T],\pi^*\Lambda^3)$ with respect to the inner product
\ben
\langle\dot\Omega_1,\dot\Omega_2\rangle_{\mb H}=\int_0^Te^{-2t}\langle\dot\Omega_1,\dot\Omega_2\rangle_{W^{1,2}}dt+\int_0^Te^{-2t}\langle\partial_t\dot\Omega_1,\partial_t\dot\Omega_2\rangle_{L^2}dt.
\ee
Note that the quadratic form $\dot \Omega \mapsto \langle - L_{\bar \Omega} \dot \Omega , \dot \Omega \rangle_{L^2}$ is defined on $W^{1,2}$ in view of equation \eqref{quad_form_L}, so that we can say that $\dot \Omega \in \mb H$ satisfies $P \dot \Omega = \dot\Phi$ for $\dot\Phi \in V^0[0,T]$ {\em weakly}, if
$$
\langle \partial_t \dot \Omega, \dot \Psi \rangle_{V^0[0,T]}  + \int_0^T e^{-2t} \langle -L_{\bar\Omega}\dot \Omega, \dot \Psi \rangle_{L^2} dt = \langle \dot \Phi, \dot\Psi \rangle_{V^0[0,T]}
$$
holds for all $\dot \Psi \in C_0^\infty(M \times (0,T),\pi^*\Lambda^3)$, the space of smooth sections of $\pi^*\Lambda^3$ which vanish near the boundary $M\times\{0,T\}$.  

\medskip

We first show that for $\dot\Phi\in V^0[0,T]$ there exists a unique weak solution in $\mb H$ to the equation 
\be\label{weaksol}
D_{\bar\Omega}F(\dot\Omega)=(0,\dot\Phi).
\ee 
Following~\cite{hp99}, we use a refined version of the
Lax--Milgram lemma, cf.~Lemma 7.8 and Theorem 7.9 in~\cite{hp99} or Theorem 16 in Chapter
10 of~\cite{fr64}. As explained in \cite{hp99}, the main point is to check {\em coercivity} of the bilinear form
\ben
A(\dot\Omega_1,\dot\Omega_2)=\int_0^Te^{-2t}\langle\partial_t\dot\Omega_1,\partial_t\dot\Omega_2\rangle_{L^2}dt
+\int_0^Te^{-2t}\langle -L_{\bar\Omega}\dot\Omega_1,\partial_t\dot\Omega_2\rangle_{L^2}dt
\ee
on $C_0^\infty(M \times (0,T),\pi^*\Lambda^3) \subset \mb H$, i.e.~to establish an estimate of the form $A(\dot\Omega,\dot\Omega)\geq C \|\dot\Omega\|^2_{\mb H}$ for some positive constant $C$. This will be a consequence of the G\r arding inequality for the operator $-L_{\bar\Omega}$ of Corollary~\ref{coercivity}. First note that $\dot\Omega$ is a solution to $P\dot\Omega=\dot\Phi$ \iff $e^{-t}\dot\Omega$ is a solution to
$(P+1)e^{-t}\dot\Omega=e^{-t}\dot\Phi$. Hence by replacing $-L_{\bar \Omega}$
by $-L_{\bar \Omega}+1$ we may assume that $-L_{\bar\Omega}$ satisfies a
strict G\r arding inequality of the form
$\langle-L_{\bar\Omega}\dot\Omega,\dot\Omega\rangle_{L^2}\geq
C\|\dot\Omega\|^2_{W^{1,2}}$ with $C$ as in
Corollary~\ref{coercivity}. Since  $\int_0^T\partial_t(e^{-2t}\langle-L_{\bar\Omega}\dot\Omega,\dot\Omega\rangle_{L^2})dt=0$ by assumption, we obtain
\ben
\int_0^Te^{-2t}\langle-L_{\bar\Omega}\dot\Omega,\partial_t\dot\Omega\rangle_{L^2}dt=\int_0^Te^{-2t}\langle
-L_{\bar\Omega}\dot\Omega,\dot\Omega\rangle_{L^2}dt\geq C\int_0^Te^{-2t}\|\dot\Omega\|^2_{W^{1,2}}dt,
\ee
whence
\ben A(\dot\Omega,\dot\Omega)\geq
C\int_0^Te^{-2t}\|\dot\Omega\|^2_{W^{1,2}}dt+\int_0^Te^{-2t}\|\partial_t\dot\Omega\|^2_{L^2}dt\geq
C\|\dot\Omega\|_{\mb H}^2,
\ee
establishing coercivity and the existence of a weak solution.

\medskip

In order to improve regularity of this weak solution one needs the following estimate:
\begin{lem}[Huisken--Polden]\label{keyestimate}
Let $s\geq0$. If $\dot\Omega\in\mb H$ is a weak solution to the equation $D_{\bar\Omega}F(\dot\Omega)=(\dot\Omega_0,\dot\Phi)$ with $\dot\Omega_0\in W^{2s+1,2}$ and $\dot\Phi\in V^s[0,T]$, then $\dot \Omega \in V^{s+1}[0,T]$ and there exists a constant $C=C(\bar\Omega,s)>0$ \st
\ben
\|\dot\Omega\|^2_{V^{s+1}[0,T]}\leq C\big(\|\dot\Omega_0\|_{W^{2s+1,2}}^2+\|\dot\Phi\|_{V^s[0,T]}^2\big).
\ee
\end{lem}
\begin{prf}
See Lemma 7.13 in~\cite{hp99}.
\end{prf}

For later use we state and prove the ensuing interior estimate:

\begin{cor}[Interior estimate]\label{interiorestimate}
Let $s\geq0$. For all $\delta>0$ there exists a constant $C=C(\delta,\bar\Omega,s)>0$ such that for $\dot\Omega\in V^{s+1}[0,T]$ one has
\ben\label{parabolicreg}
\|\dot\Omega\|^2_{V^{s+1}[\delta,T]}\leq C\big(\|\dot\Omega\|^2_{V^0[0,T]} + \|P\dot\Omega\|^2_{V^s[0,T]}\big).
\ee
\end{cor}
\begin{prf}
Since $\dot\Phi=P\dot\Omega$ is of class $V^s$, Lemma~\ref{keyestimate} gives 
\ben
\|\dot\Omega\|^2_{V^{s+1}[0,T]}\leq
C(\bar\Omega,s)\bigl(\|\dot\Omega_0\|^2_{W^{2s+1,2}}+\|P\dot\Omega\|^2_{V^s[0,T]}\bigr). 
\ee
Put $\Omega':=\varphi\cdot\dot\Omega$, where $\varphi:[0,T]\to\R$ is a smooth cut--off function satisfying $\varphi(t)=0$ for $t$ close to $0$ and $\varphi(t)=1$ for $t\in[\epsilon,T]$. Then $P\Omega'=\partial_t\varphi\cdot\dot\Omega+\varphi\cdot P\dot\Omega$ and $\Omega'_0=0$, whence
\begin{eqnarray*}
\|\dot\Omega\|_{V^{s+1}[\epsilon,T]} & \leq &\|\Omega'\|_{V^{s+1}[0,T]}\leq C(\bar\Omega,s)^{\frac{1}{2}} \|\partial_t\varphi\cdot\dot\Omega+\varphi\cdot P\dot\Omega\|_{V^s[0,T]}\\
& \leq & C(\bar\Omega,s)^{\frac{1}{2}}\bigl({\textstyle\sup_{t\in[0,T]}|\partial_t\varphi(t)|}\cdot\|\dot\Omega\|_{V^s[0,T]}+\|P\dot\Omega \|_{V^s[0,T]}\bigr).  
\end{eqnarray*}
Induction on $s$ yields the result.
\end{prf}

\medskip

Mutatis mutandis one proves along the lines of Theorem 7.14 in~\cite{hp99}:

\begin{thm}[Huisken--Polden]\label{banach}
The map $D_{\bar\Omega}F$ in~\eqref{banachflowlin} is a Banach space isomorphism.
\end{thm}

As a consequence we obtain existence of the Dirichlet--DeTurck flow together
with the a priori estimate on a finite time interval $[0,T]$ for initial conditions
(depending a priori on $T$) sufficiently close to $\bar \Omega$:

\begin{cor}\label{T-existence}
For all $\epsilon>0$ and $0<T<\infty$, there exists $\delta=\delta(\bar\Omega,\epsilon,T)>0$
such that for $\widetilde\Omega_0\in W^{k,2}$ with
$\|\widetilde\Omega_0-\bar\Omega\|_{W^{k,2}}<\delta$, the Dirichlet--DeTurck
flow $\widetilde\Omega_t$ at $\bar\Omega$ with initial condition $\widetilde\Omega_0$ exists and satisfies $\|\widetilde\Omega_t-\bar\Omega\|_{W^{k,2}}<\epsilon$ for all $t\in[0,T]$.
\end{cor}
\begin{prf}
Since $k$ is odd, $s := (k+1)/2$ is an integer. Let $\bar\Omega$ also denote the pull--back to $M\times[0,T]$. Then $\bar\Omega\in V^s[0,T]$ and $F(\bar\Omega)=(\bar\Omega,0)$. By virtue of the previous Lemma~\ref{banach} and the Banach space inverse function theorem, the map $F$ in~\eqref{banachflowmap} is a local diffeomorphism near $\bar\Omega$. The trace maps $\Omega\in V^s[0,T]\mapsto \Omega_t\in W^{k,2}$ are continuous with a uniform bound on their norms since $t$ varies within a compact interval. Hence there exists $C>0$ such that $\|\widetilde\Omega_t-\bar\Omega\|_{W^{k,2}} \leq C \|\widetilde\Omega-\bar\Omega\|_{V^s[0,T]}$ holds for all $t\in[0,T]$. For suitably chosen $\delta>0$, the condition $\|\widetilde\Omega_0-\bar\Omega\|_{W^{k,2}}<\delta$ implies that $(\widetilde\Omega_0,0)$ is close enough to $(\bar\Omega,0)$ to ensure that $\widetilde\Omega=F^{-1}(\widetilde\Omega_0,0)$ satisfies $\|\widetilde\Omega-\bar\Omega\|_{V^s[0,T]} < \epsilon/C$. 
\end{prf}

\begin{rem}
Using Theorem \ref{banach} it is also possible to give an alternative proof of short--time existence for the Dirichlet--DeTurck flow, cf.~\cite{hp99} for details.
\end{rem}

\medskip

As in Section~\ref{linstab} let $\mc{S}_{\bar\Omega}\subset\widetilde Q^{-1}_{\bar \Omega}(0)$
denote a suitably chosen slice around $\bar \Omega$. For a positive $3$--form $\bar\Omega'$ close to $\bar \Omega$ we write $\widetilde \Omega_t = \bar \Omega'+\omega_t'$. Let now 
\ben
L_{\bar {\Omega}'}:=D_{\bar {\Omega}'}\widetilde{Q}_{\bar {\Omega}}
\ee
and
\ben
R_{\bar \Omega'}(\omega_t'):=\widetilde Q_{\bar \Omega}(\widetilde\Omega_t)-L_{\bar \Omega'}\omega_t'.
\ee
Then we can recast the flow equation into
\be\label{flowequation}
\frac{\partial}{\partial t}\widetilde{\Omega}_t =  L_{\bar {\Omega}'}
\omega_t' + R_{\bar {\Omega}'}(\omega_t').
\ee
The basic idea is that the behaviour of the Dirichlet--DeTurck flow should be
dominated by the linear term $L_{\bar \Omega'}$. To obtain precise results, we
need to control the remainder term $R_{\bar\Omega'}$.

\bigskip

In order to analyse $R_{\bar {\Omega}'}$ we introduce the following notation: Let $E_1, E_2, E$ and $F$ be tensor bundles over $M$ equipped with the bundle metrics induced by $\bar \Omega$. Denote by $\bar {\nabla}$ the covariant derivative associated with $\bar
{\Omega}$. Let $\Omega$ be a further positive $3$--form of class $C^2$ and let $\bar\Omega' \in
\mc{S}_{\bar\Omega}$. Put $\omega'=\Omega-\bar\Omega'$ and assume that
$\|\omega'\|_{C^2}<\epsilon_1$ and $\|\bar\Omega'-\bar\Omega\|_{C^2}<\epsilon_2$. For $s\in\Gamma(M,E)$ we will generically write $\circledast_l s$ for sections of the form $A(s)\in\Gamma(M,F)$, where $A\in\Gamma(M,E^*\otimes F)$ is a section of class $C^l$ possibly depending on $\omega'$ and $\bar\Omega'$ \st $\|A\|_{C^l}\leq C(\epsilon_1,\epsilon_2)$. Using this convention, we have for instance 
\ben
\circledast_k\circledast_ls=\circledast_ks
\ee
for $k\leq l$ and, by the product rule for $\bar\nabla$,
\ben
\bar\nabla(\circledast_ls)=\circledast_{l-1}s+\circledast_l\bar\nabla s.
\ee
Similarly, for a section
\ben
B \in \Gamma (M,E_1^* \otimes E_2^* \otimes F)
\ee
with $\|B\|_{C^l}\leq C(\epsilon_1,\epsilon_2)$ and $s_i \in \Gamma(M,E_i)$ we write generically
$s_1 \circledast_l s_2 $ for the section $B(s_1,s_2) \in \Gamma(M,F)$. For instance we have
\ben
(\circledast_{l_0}s_0)\circledast_k(\circledast_{l_1}s_1)=s_1\circledast_ks_2,\quad\mbox{and}\quad\circledast_k(s_1\circledast_ls_2)=s_1\circledast_ks_2
\ee
for $k\leq l,\,l_0,\,l_1$. In all cases we set $\circledast:=\circledast_0$. As we will differentiate at most twice we take $l$ to be smaller or equal than two even if we have $C^l$--control for higher $l$.

\begin{ex}
Since the exterior differential $d$ is obtained by anti--symmetrising the covariant derivative $\bar \nabla$, we have $d\omega'=\circledast_2\bar\nabla\omega'$. Hence $\star_{\bar \Omega'}d\omega' =\circledast_2\bar\nabla\omega'$ and we get
\begin{eqnarray*}
d \star_{\bar \Omega'} d\omega' &=&\circledast_2\bar\nabla(\circledast_2\bar\nabla\omega')\\
&=&\circledast_2(\circledast_1 \bar \nabla \omega' + \circledast_2 \bar\nabla^2
\omega')\\
&=&\circledast_1 \bar \nabla \omega' + \circledast_2 \bar\nabla^2
\omega'.
\end{eqnarray*}
\end{ex}

The following result gives a rough description of the structure of the
remainder term, which however will turn out to be sufficient for our
purposes.

\begin{lem}\label{symb_calc}
Let $\Omega$ be a positive $3$--form of class $C^2$ and $\bar\Omega' \in
\mc{S}_{\bar\Omega}$. Let further $\omega'=\Omega-\bar\Omega'$. Assume that
$\|\omega'\|_{C^2}<\epsilon_1$ and
$\|\bar\Omega'-\bar\Omega\|_{C^2}<\epsilon_2$. Then 
\be\label{symbcalcforR}
R_{\bar\Omega'}(\omega') =  \omega' \circledast \omega' + \omega' \circledast \bar\nabla \omega' +\omega' \circledast
\bar\nabla^2 \omega'+ \bar\nabla \omega' \circledast \bar \nabla \omega'.
\ee
In particular, there exists $C=C(\epsilon_1,\epsilon_2)>0$ such that 
\ben
|R_{\bar \Omega'}(\omega')| \leq \epsilon_1C ( |\omega'| + |\bar {\nabla} \omega'|).
\ee
\end{lem}

\begin{rem}
In particular we absorbed into $\circledast$ any term in~\eqref{symbcalcforR} of order strictly higher than two in $\omega'$.
\end{rem}

\begin{prf}
Recall that $\widetilde{Q}_{\bar {\Omega}}(\Omega)=Q(\Omega)+\Lambda_{\bar {\Omega}}(\Omega)$
where
\ben
Q(\Omega) = -\big(\Delta_{\Omega}\Omega+\frac{1}{3}[d\delta_{\Omega}\Omega]_1-2[d\delta_{\Omega}\Omega]_{27}+q_{\bar \Omega}(\bar \nabla \Omega)\big),
\ee
projections being taken with respect to $\bar\Omega$, and
\ben
\Lambda_{\bar {\Omega}}(\Omega)=\mathcal{L}_{X(\Omega)}\Omega
\ee
with $X(\Omega)=-(\delta_{\bar\Omega}\Omega)\llcorner\bar \Omega$. In the
following we will calculate the difference between $\widetilde Q_{\bar\Omega}$ and its linearisation at $\bar \Omega'$ term by term. First write $\Omega=\bar\Omega'+\omega'$ and observe that mapping a positive 3--form $\Omega$ to the corresponding Hodge operator $\star_\Omega$ gives rise to a fibre--preserving smooth map $\star: \Lambda^3_+T^*M
\rightarrow \End(\Lambda^*T^*M)$. Hence there exists a fibrewise linear map  $A_{\bar \Omega'}(\omega'):
\Lambda^3T^*M \rightarrow \End(\Lambda^*T^*M)$, depending smoothly on
$\bar\Omega'$ and $\omega'$, such that
$\star_\Omega=\star_{\bar\Omega'}+A_{\bar \Omega'}(\omega')\omega'$; in
particular $(D_{\bar\Omega'}\star)(\omega')=A_{\bar\Omega'}(0)\omega'$.
Using $d\bar\Omega'=0$, we get
\begin{eqnarray*}
\delta_{\Omega}d \Omega &=& \star_{\Omega}d\star_{\Omega}d\omega'\\
&=& \star_{\bar\Omega'} d\star_{\bar\Omega'}d \omega' + \star_{\bar\Omega'}d (A_{\bar \Omega'}(\omega')\omega') d\omega' + (A_{\bar \Omega'}(\omega')\omega') d\star_{\bar\Omega'}d\omega'\\
&&+ (A_{\bar \Omega'}(\omega')\omega') d (A_{\bar \Omega'}(\omega')\omega') d\omega'.
\end{eqnarray*}
Subtracting the linear term in $\omega'$, we get 
\begin{eqnarray*}
& &\delta_{\Omega}d \Omega-\star_{\bar\Omega'} d\star_{\bar \Omega'}d \omega'\\
&=&\star_{\bar \Omega'}d(A_{\bar \Omega'}(\omega')\omega') d\omega' + (A_{\bar \Omega'}(\omega')\omega')d\star_{\bar\Omega'}d\omega' + (A_{\bar \Omega'}(\omega')\omega') d (A_{\bar \Omega'}(\omega')\omega')d\omega'\\
&=&\omega' \circledast \bar\nabla \omega' +\omega' \circledast
\bar\nabla^2 \omega'+ \bar\nabla \omega' \circledast \bar \nabla \omega'. 
\end{eqnarray*}
Similarly, using $\delta_{\bar\Omega'}\bar\Omega'=0$, we get
\begin{eqnarray*}
d \delta_{\Omega}\Omega &=&-d\star_{\Omega}d\star_{\Omega}\Omega\\
&=& - d\star_{\bar \Omega'} d\star_{\bar \Omega'} \omega' - d \star_{\bar \Omega'}d (A_{\bar \Omega'}(\omega')\omega')\bar\Omega'-d\star_{\bar\Omega'} d (A_{\bar \Omega'}(\omega')\omega') \omega' \\
& & - d (A_{\bar \Omega'}(\omega')\omega')d \star_{\bar \Omega '} \omega'  -d  (A_{\bar \Omega'}(\omega')\omega') d (A_{\bar\Omega'}(\omega')\omega') \bar\Omega'\\
& & - d  (A_{\bar \Omega'}(\omega')\omega')d(A_{\bar \Omega'}(\omega')\omega')\omega'. 
\end{eqnarray*}
Again, subtracting the linear term in $\omega'$, we get 
\begin{eqnarray*}
& &d\delta_{\Omega}\Omega+d\star_{\bar\Omega'}d\star_{\bar\Omega'}\omega'+d\star_{\bar\Omega'}d(A_{\bar \Omega'}(0)\omega')\bar\Omega'\\
&=&-d \star_{\bar \Omega'}d((A_{\bar \Omega'}(\omega')-A_{\bar\Omega'}(0))\omega')\bar\Omega' - d \star_{\bar\Omega'} d (A_{\bar \Omega'}(\omega')\omega') \omega'\\
& & - d (A_{\bar \Omega'}(\omega')\omega')d \star_{\bar \Omega '} \omega'  - d(A_{\bar \Omega'}(\omega')\omega')  d (A_{\bar \Omega'}(\omega')\omega') \bar\Omega'\\
& & - d (A_{\bar \Omega'}(\omega')\omega') d (A_{\bar\Omega'}(\omega')\omega') \omega' \\
&=& \omega' \circledast \omega' + \omega' \circledast \bar\nabla \omega' +\omega' \circledast
\bar\nabla^2 \omega'+ \bar\nabla \omega' \circledast \bar \nabla \omega'.  
\end{eqnarray*}
This takes care of the first three terms in $Q$ (note that the linear
projections onto irreducible components are absorbed into $\circledast$). Furthermore, the term coming from the quadratic form $q_{\bar\Omega}$ contributes a remainder term of type $\bar \nabla \omega' \circledast \bar\nabla \omega'$. Finally, in order to deal with the term $\Lambda_{\bar\Omega}(\Omega)$, we observe that $X(\bar \Omega')=0$ since $\bar\Omega'\in\mc{S}_{\bar\Omega}$ which implies $[\delta_{\bar\Omega}\bar\Omega']_7$=0. Hence $X(\Omega)=X(\omega')$
and
\ben
\mathcal{L}_{X(\omega')}\omega' =\omega'\circledast \bar\nabla \omega' + \omega'\circledast \bar\nabla^2 \omega' + \bar
\nabla \omega' \circledast \bar \nabla \omega'.
\ee
This finishes the proof.
\end{prf}

In the following we need some standard results from perturbation theory of linear
operators. This is summarised in the following statement:

\begin{lem}\label{perturb}
For all $\epsilon>0$ there exists $\delta=\delta(\epsilon)>0$ such that if $\|\bar \Omega' - \bar {\Omega}\|_{W^{k,2}}<\delta$, then
\ben
\langle -L_{\bar \Omega'}\omega,\omega\rangle_{L^2} \geq(1-\epsilon)\langle-L_{\bar\Omega}\omega,\omega\rangle_{L^2}-\epsilon\|\omega\|^2_{L^2}
\ee
for all $\omega \in W^{2,2}$.
\end{lem}
\begin{prf}
We apply Theorem 9.1 in \cite{we80} to the operators
$T:=-\epsilon L_{\bar \Omega}$ and $V:=L_{\bar\Omega}-L_{\bar\Omega'}$. Without loss of generality we may assume that $L_{\bar\Omega'}$ is symmetric; otherwise we replace it by the symmetric
operator $\frac{1}{2}(L_{\bar\Omega'}+L_{\bar\Omega'}^*)$, where
$*$ denotes the formal adjoint taken with respect to $\bar\Omega$.
Then by elliptic regularity for the operator $T$ the estimate
\ben
\|V\omega\|_{L^2} \leq a \|\omega \|_{L^2} + b \|T\omega\|_{L^2}
\ee
holds for arbitrarily small $a,b >0$ if $\bar \Omega'$ is sufficiently close to $\bar \Omega$ \wrt the $W^{k,2}$--norm. In particular, $V$ will be $T$--bounded with $T$--bound less than $1$. Since $T$ is non--negative we obtain with $T+V=-L_{\bar\Omega'}+(1-\epsilon)L_{\bar \Omega}$ that
\ben
-L_{\bar \Omega'}+(1-\epsilon)L_{\bar \Omega} \geq -\epsilon
\ee
or equivalently
\ben
\langle -L_{\bar{\Omega}'}\omega,\omega\rangle_{L^2} \geq(1-\epsilon)\langle-L_{\bar\Omega}\omega,\omega\rangle_{L^2}-\epsilon\|\omega\|^2_{L^2}
\ee
for all $\omega \in W^{2,2}$ and $\bar \Omega'$ sufficiently close to
$\bar \Omega$.
\end{prf}

The family $(L_{\bar\Omega'})_{\bar\Omega'\in\mc{S}_{\bar\Omega}}$ is a smooth family of
elliptic operators, and hence gives rise to a smooth family of Fredholm
operators $L_{\bar\Omega'}:W^{2,2} \rightarrow L^2$. Since
$\mc{S}_{\bar\Omega}$ is a smooth (finite--dimensional) manifold and
$\mc{S}_{\bar\Omega} \subset \widetilde{Q}_{\bar\Omega}^{-1}(0)$, we have
$T_{\bar {\Omega}'}\mc{S}_{\bar\Omega} \subset \ker L_{\bar {\Omega}'}$. Furthermore $T_{\bar {\Omega}}\mc{S}_{\bar\Omega} = \ker L_{\bar {\Omega}}$ by Proposition \ref{secvar}. Usual Fredholm theory implies that $T_{\bar {\Omega}'}\mc{S}_{\bar\Omega} = \ker L_{\bar {\Omega}'}$ for $\bar\Omega'$ sufficiently close to $\bar\Omega$, and hence the kernels and ranges of $L_{\bar\Omega'}$ form smooth vector bundles (the latter infinite--dimensional). Interpreting $(\ker L_{\bar\Omega'})^\perp$ as the fibre at $\bar\Omega'$ of the normal bundle of $\mc{S}_{\bar\Omega}$ with respect to the $L^2$--metric, we obtain the following statement:

\begin{lem}[Orthogonal projection]\label{orth_proj}
There exists $\epsilon>0$ such that if $\|\Omega-\bar\Omega\|_{W^{k,2}} < \epsilon$, then there exists a unique $\bar \Omega' \in \mc{S}_{\bar \Omega}$ such that $\omega'=\Omega - \bar\Omega' \in (\ker L_{\bar \Omega'})^\perp$, the orthogonal complement being taken
inside $L^2$.
\end{lem}

\begin{rem}
In the situation of Lemma \ref{orth_proj} we obtain using the continuity of the embedding $W^{k,2} \hookrightarrow C^2$ that there exist constants $\epsilon_1,\epsilon_2>0$ such that $\|\Omega-\bar\Omega\|_{W^{k,2}}<\epsilon$ implies that $\|\omega'\|_{C^2} < \epsilon_1$ and $\|\bar\Omega'-\bar\Omega\|_{C^2}<\epsilon_2$ with $\epsilon_1,\epsilon_2 \rightarrow 0$ as $\epsilon \rightarrow 0$. 
\end{rem}

\begin{prp}\label{estimate_range}
For all $\kappa>0$ there exists $\epsilon>0$ such that if $\|\Omega-\bar\Omega\|_{W^{k,2}} < \epsilon$, then
\ben
\|R_{\bar\Omega'}(\omega')\|_{L^2} \leq \kappa \|L_{\bar \Omega'}\omega'\|_{L^2}
\ee
with $\bar\Omega'\in\mc{S}_{\bar\Omega}$ and $\omega'=\Omega - \bar\Omega'$ as in Lemma \ref{orth_proj}.
\end{prp}
\begin{prf}
Elliptic regularity implies that
\ben
L_{\bar \Omega'} : (\ker L_{\bar \Omega'})^\perp \cap W^{2,2}\rightarrow \operatorname{im} L_{\bar \Omega'} \subset L^2
\ee
is a Banach space isomorphism, the orthogonal complement being taken inside $L^2$. In particular, there exists $C>0$ such that $\|\omega\|_{W^{2,2}} \leq C \|L_{\bar \Omega'}\omega\|_{L^2}$ for all $\omega \in (\ker L_{\bar \Omega'})^\perp \cap W^{2,2}$. This constant $C$ can be chosen uniform in $\bar \Omega'$, since $L_{\bar\Omega'}$ is a smooth family of elliptic operators parametrised by a finite--dimensional, hence locally compact manifold. Using Lemma \ref{symb_calc} we get
\ben
\|R_{\bar \Omega'}(\omega')\|_{L^2} \leq \epsilon_1C'
\|\omega'\|_{W^{2,2}}\leq \epsilon_1CC'  \|L_{\bar \Omega}\omega'\|_{L^2}
\ee
with $\epsilon_1,\epsilon_2>0$ as in the remark following Lemma \ref{orth_proj} and some
constant $C'=C'(\epsilon_1,\epsilon_2)>0$. For $\epsilon_1\ll 1$, this constant $C'$ can in fact be chosen independent of $\epsilon_1$ in the above
estimate, which implies the result.
\end{prf}

As a consequence, we obtain that locally $\widetilde Q_{\bar\Omega}^{-1}(0)$ coincides with $\mc{S}_{\bar\Omega}$, more precisely:

\begin{cor}\label{zero_set}
There exists $\epsilon>0$ such that if $\|\Omega-\bar\Omega\|_{W^{k,2}} < \epsilon$, then $\widetilde{Q}_{\bar\Omega}(\Omega)=0$ implies that $\Omega \in \mc{S}_{\bar\Omega}$.
\end{cor}
\begin{prf}
Let $\bar\Omega' \in \mc{S}_{\bar\Omega}$ be chosen according to Lemma
\ref{orth_proj}, i.e.~such that $\omega'=\Omega-\bar\Omega'$ satisfies $\omega' \in (\ker L_{\bar\Omega'})^\perp$. We write $\widetilde Q_{\bar \Omega}(\Omega) = L_{\bar \Omega'} \omega' + R_{\bar\Omega'}(\omega')$. Now $\widetilde Q_{\bar \Omega}(\Omega)=0$ implies that $\|R_{\bar\Omega'}(\omega')\|_{L^2} = \|L_{\bar \Omega'}\omega'\|_{L^2}$, which in view of Proposition \ref{estimate_range} (with e.g.~$\kappa=\frac{1}{2}$) yields $R_{\bar\Omega'}(\omega')=L_{\bar \Omega'}\omega'=0$. Since $\omega'$ was chosen orthogonal to the kernel of $L_{\bar\Omega'}$, we obtain $\omega'=0$, hence $\Omega \in \mc{S}_{\bar\Omega}$.
\end{prf}

Let $\widetilde \Omega_t$ be a Dirichlet--DeTurck flow solution on $[0,T]$ and let
$\widetilde Q_t:=\widetilde Q_{\bar \Omega}(\widetilde \Omega_t)$. Then
$\frac{\partial}{\partial t} \widetilde \Omega_t = \widetilde Q_t$ and
$\widetilde Q_t$ satisfies the linearised flow equation
\begin{equation}\label{linearised_flow_eqnI}
\frac{\partial}{\partial t} \widetilde Q_t = L_{\widetilde\Omega_t} \widetilde Q_t
\end{equation}
with $L_{\widetilde\Omega_t}=D_{\widetilde \Omega_t}\widetilde Q_{\bar\Omega}$.
This is a linear parabolic equation with time--dependent
coefficients. We view the operator $L_{\widetilde\Omega_t}$ as a
(non--symmetric) perturbation of the symmetric operator $L_{\bar \Omega}$, to which in particular Lemma \ref{perturb} applies.

\begin{lem}[$L^2$--almost orthogonality of $\widetilde Q_t$]\label{almost_orth} 
For all $\kappa>0$ there exists $\epsilon > 0$ such that if $\|\widetilde\Omega_t-\bar
\Omega\|_{W^{k,2}}<\epsilon$ for all $t \in [0,T]$, then
\ben
|\, \langle \widetilde Q_t, \omega \rangle_{L^2}| \leq \kappa \,\|
\widetilde Q_t\|_{L^2}\|\omega\|_{L^2}
\ee
for all $\omega \in \ker L_{\bar \Omega}$.
\end{lem}
\begin{prf}
Let $\bar\Omega_t \in \mc{S}_{\bar\Omega}$ be chosen according to Lemma
\ref{orth_proj}, i.e.~such that $\omega_t'=\widetilde\Omega_t-\bar\Omega_t$ satisfies $\omega_t' \in (\ker L_{\bar\Omega_t})^\perp$. We write
$\widetilde Q_t = L_{\bar \Omega_t} \omega_t' +
R_{\bar\Omega_t}(\omega_t')$. First we observe that due to the symmetry of $L_{\bar \Omega}$ its range is $L^2$--orthogonal to its
kernel. Since the ranges of the operators $L_{\bar\Omega'}$ form a smooth vector bundle for $\bar\Omega' \in \mc{S}_{\bar\Omega}$, this implies that they are $L^2$--almost orthogonal to $\ker L_{\bar \Omega}$ in the above sense, in particular
\ben
|\, \langle L_{\bar\Omega_t}\omega_t', \omega \rangle_{L^2}| \leq \kappa \,\|
 L_{\bar\Omega_t}\omega_t'\|_{L^2}\|\omega\|_{L^2}
\ee
for all $\omega \in \ker L_{\bar\Omega}$, if $\epsilon>0$ is chosen sufficiently small. Proposition \ref{estimate_range} then shows that
$$
|\, \langle R_{\bar\Omega_t}(\omega_t'), \omega \rangle_{L^2}| \leq \| R_{\bar\Omega_t}(\omega_t') \|_{L^2} \|\omega\|_{L^2} \leq \kappa \|  L_{\bar\Omega_t}\omega_t' \|_{L^2} \|\omega\|_{L^2} 
$$
and hence
\begin{align*}
|\, \langle \widetilde Q_t, \omega \rangle_{L^2}| &\leq |\, \langle L_{\bar\Omega_t}\omega_t', \omega \rangle_{L^2}| + |\, \langle R_{\bar\Omega_t}(\omega_t'), \omega \rangle_{L^2}|\\
& \leq  2\kappa \,\| L_{\bar\Omega_t}\omega_t' \|_{L^2}\|\omega\|_{L^2}\leq 2\kappa(1-\kappa)^{-1}\|
\widetilde Q_t\|_{L^2}\|\omega\|_{L^2}. 
\end{align*}
for all $\omega \in \ker L_{\bar\Omega}$, again if $\epsilon>0$ is chosen sufficiently small. This proves the result.
\end{prf}

\begin{lem}[$L^2$--exponential decay of $\widetilde{Q}_t$]\label{L2_decay}
There exist $\epsilon>0$ and $\lambda=\lambda(\bar\Omega)>0$ such that if $\|\widetilde\Omega_t-\bar
\Omega\|_{W^{k,2}}<\epsilon$ for all $t \in [0,T]$, then
\ben
\|\widetilde Q_t\|^2_{L^2}\leq e^{-\lambda t}\|\widetilde Q_0\|^2_{L^2}.
\ee
for all $t \in [0,T]$. 
\end{lem}
\begin{prf}
Using equation~\eqref{linearised_flow_eqnI} we obtain 
\ben
\frac{d}{dt} \frac{1}{2} \|\widetilde Q_t\|^2_{L^2} = \langle
L_{\widetilde \Omega_t}\widetilde Q_t , \widetilde Q_t\rangle_{L^2}
\ee
Let $\lambda_1 >0$ be the first (positive) eigenvalue of $-L_{\bar {\Omega}}$. Then one has
$\langle -L_{\bar {\Omega}} \omega, \omega \rangle_{L^2} \geq \lambda_1\|\omega\|^2_{L^2}$
for all $\omega \in (\ker L_{\bar {\Omega}})^\perp$, and hence, using Lemma
\ref{almost_orth}, that
\ben
\langle -L_{\bar {\Omega}}\widetilde Q_t , \widetilde Q_t \rangle_{L^2} \geq \frac{3\lambda_1}{4}
\|\widetilde Q_t\|^2_{L^2}
\ee
for $\epsilon>0$ sufficiently small. Using Lemma \ref{perturb} one gets
\ben
\langle -L_{\widetilde\Omega_t} \widetilde Q_t, \widetilde Q_t \rangle_{L^2} \geq
(1-\epsilon') \frac{3\lambda_1}{4} \|\widetilde Q_t\|^2_{L^2} - \epsilon'
\|\widetilde Q_t\|^2_{L^2}\geq \frac{\lambda_1}{2} \|\widetilde Q_t\|^2_{L^2}
\ee
if $\epsilon,\epsilon'>0$ are chosen sufficiently small. Finally we get
\ben
\frac{d}{dt} \frac{1}{2} \|\widetilde Q_t\|^2_{L^2} \leq - \lambda \|\widetilde Q_t\|^2_{L^2}
\ee
with $\lambda:=\lambda_1/2$ and $\epsilon>0$ sufficiently small. Now Gronwall's lemma implies the result.
 \end{prf}

Parabolic regularity theory now yields higher derivative estimates:

\begin{prp}[$W^{k,2}$--exponential decay of $\widetilde{Q}_t$]\label{Wk2_decay}
There exist $\epsilon>0$ and $C,\lambda>0$ such that if $\|\widetilde\Omega_t-\bar
\Omega\|_{W^{k,2}}<\epsilon$ for all $t \in [0,T]$, then
\ben
\|\widetilde Q_t\|^2_{W^{k,2}}\leq Ce^{-\lambda t}
\ee
for all $t \in [0,T]$. 
\end{prp}
\begin{prf}
If $\epsilon>0$ and $\lambda >0$ are chosen according to Lemma \ref{L2_decay}, then for $0 \leq t_0 < t_1 \leq T$ we get
\ben
\int_{t_0}^{t_1} \|\widetilde Q_\tau\|_{L^2}^2 d\tau\leq \|\widetilde Q_0\|_{L^2}^2\int_{t_0}^{t_1} e^{-\lambda \tau}d\tau \leq C_0e^{-\lambda t_0}
\ee
for $C_0:=\|\widetilde Q_0\|_{L^2}^2 /\lambda$. It is easy to see that the estimate of Corollary \ref{interiorestimate} remains true for $P_t:= \partial_t-L_{\widetilde\Omega_t}$ (with a constant $C_s=C_s(\epsilon,\delta,\bar \Omega)>0$) if $\epsilon>0$ is sufficiently small. In particular we get for $s:=(k-1)/2$ that
\ben
\|\widetilde Q\|^2_{V^{s+1}[t_0+\delta,t_1]}\leq C_s \|\widetilde Q\|^2_{V^0[t_0,t_1]}\leq C_s \int_{t_0}^{t_1}\|\widetilde Q_\tau\|_{L^2}^2 d\tau 
\ee
since $P_t\widetilde Q_t=0$. Combining these two estimates and using the trace theorem as in the proof of Corollary \ref{T-existence} yields the result.
\end{prf}

We will now describe the choice of the constants:

\begin{enumerate}
	\item Choose $\epsilon > 0$ such that
	\begin{enumerate}
		\item $\|\widetilde\Omega_0-\bar \Omega\|_{W^{k,2}} < \epsilon$ implies that the Dirichlet--DeTurck flow at $\bar\Omega$ with initial condition $\widetilde{\Omega}_0$ exists on $[0,1]$. This is possible using Lemma \ref{T-existence}.
		\item $\|\widetilde\Omega_t-\bar \Omega\|_{W^{k,2}} < \epsilon$ for $t \in [0,T]$ yields
\ben
\|\widetilde Q_t\|_{W^{k,2}} \leq Ce^{-\lambda t}
\ee 
on $[0,T]$ for any $T<T_{max}$, where $T_{max}$ is the maximal life--time of the flow. This is possible using Proposition \ref{Wk2_decay}. 
\end{enumerate}
	\item Choose $T>0$ such that $\int_T^\infty C e^{-\lambda t} dt < \epsilon/2$ with $C$ and $\lambda$ as above.
	\item Choose $\delta >0$ such that
\ben
\|\widetilde\Omega_0-\bar \Omega\|_{W^{k,2}} < \delta
\ee
implies that the Dirichlet--DeTurck flow at $\bar\Omega$ with initial condition $\widetilde{\Omega}_0$ exists on $[0,T]$ with $\|\widetilde\Omega_t-\bar \Omega\|_{W^{k,2}}< \epsilon/2$. This again is possible using Lemma \ref{T-existence}.
\end{enumerate}

The proof of Theorem \ref{stability} may then be finished as follows:

\medskip

Let the initial condition $\widetilde\Omega_0$ be given satisfying
\ben
\|\widetilde\Omega_0-\bar \Omega\|_{W^{k,2}}< \delta.
\ee
Let $T_{max}>T$ such that $[0,T_{max})$ is the maximal time interval on which the Dirichlet--DeTurck flow with initial condition $\widetilde\Omega_0$ exists. Suppose now that $T_{max} < \infty$. For $T \leq t<T_{max}$ one has
\ben
\widetilde\Omega_t=\widetilde\Omega_T + \int_T^t \frac{\partial}{\partial \tau}\widetilde\Omega_\tau \,d\tau
\ee
and hence, as $\frac{\partial}{\partial t}\widetilde \Omega_t=\widetilde Q_t$,  
\begin{eqnarray*}
\|\widetilde\Omega_t-\widetilde\Omega_T\|_{W^{k,2}} &=& \Bigl\|\int_T^t
\widetilde Q_\tau d\tau \Bigr\|_{W^{k,2}}\\
&\leq&\int_T^t \| \widetilde Q_\tau\|_{W^{k,2}} d\tau\\
&\leq&\int_T^\infty C e^{-\lambda \tau} d\tau < \epsilon/2.
\end{eqnarray*}
Using this for $t_0=\max\{T,T_{max}-1/2\}$ and bearing the assumption $\|\widetilde\Omega_T-\bar
\Omega\|_{W^{k,2}}<\epsilon/2$ in mind, we obtain
\ben
\|\widetilde\Omega_{t_0}-\bar \Omega\|_{W^{k,2}}\leq
\|\widetilde\Omega_{t_0}-\widetilde\Omega_T\|_{W^{k,2}} +
\|\widetilde\Omega_T-\bar \Omega\|_{W^{k,2}}<\epsilon.
\ee
Therefore the Dirichlet--DeTurck flow can be continued at least up to time
$t_0+1\geq T_{\max}+1/2$, a contradiction. Hence $T_{max}=\infty$ and we have
established longtime--existence together with the estimate
$\|\widetilde\Omega_t-\bar \Omega\|_{W^{k,2}} < \epsilon$ for all $t \in [0,\infty)$.

\medskip

Finally we set
\ben
\widetilde\Omega_\infty=\widetilde\Omega_0+\int_0^\infty\frac{\partial}{\partial t}\widetilde\Omega_t\,dt.
\ee
Using Proposition \ref{Wk2_decay} we see that this integral converges,
i.e.\ that $\widetilde\Omega_t$ converges to $\widetilde\Omega_\infty$ with
respect to the $W^{k,2}$--norm. 

\medskip

Furthermore, since $\widetilde Q_{\bar \Omega}:W^{k,2} \rightarrow W^{k-2,2}$ is continuous (even differentiable), we obtain 
\ben
0=\lim_{t \rightarrow \infty}\widetilde Q_t=\widetilde Q_{\bar \Omega}(\widetilde \Omega_\infty),
\ee
in other words, $\widetilde \Omega_\infty \in \widetilde Q_{\bar \Omega}^{-1}(0)$. Since locally $\widetilde{Q}_{\bar \Omega}^{-1}(0)=\mc{S}_{\bar \Omega}$ according to Corollary~\ref{zero_set}, we obtain $\widetilde\Omega_\infty \in \mc{S}_{\bar \Omega}$. In particular, $\widetilde\Omega_\infty$ is a torsion--free $\Gt$--form. This finishes the proof of Theorem \ref{stability}.
\hfill$\blacksquare$

\bigskip

\centerline{\textbf{Acknowledgements}}

\medskip

The authors would like to thank Richard Bamler, Nata\v sa \v Se\v sum, Miles Simon and Jan Swoboda for explaining their related work to us. 
%
%
%
%
%
\appendix\section{Appendix: $\Gt$--differential operators}\label{appendix}
The $\Gt$--differential operators $d^p_q$ we use in the text are given as follows:

\ben
d^1_q:\Omega_1\to\Omega_q\quad
\left\{\begin{array}{lcl}
d^1_1f & \equiv & 0\\[2pt]
d^1_{14}f & \equiv & 0\\[2pt]
d^1_7f & = & df\\[2pt]
d^1_{27}f & \equiv & 0
\end{array}\right.
\ee

\bigskip

\ben
d^7_q:\Omega_1\to\Omega_q\quad
\left\{\begin{array}{lcl}
d^7_1\alpha & = & \delta_\Omega\alpha=\tfrac{1}{4}\star_\Omega\! d\big(\star_\Omega\!(\alpha\wedge\Omega)\wedge\Omega\big)\\[2pt]
d^7_7\alpha & = & \star_\Omega(d\alpha\wedge\star_\Omega\Omega)=-\tfrac{1}{2}\star_\Omega\!\big(\delta_\Omega(\alpha\wedge\Omega)\wedge\Omega\big)\\[2pt]
d^7_{14}\alpha & = & [d\alpha]_{14}\\[2pt]
d^7_{27}\alpha & = & [d\star_\Omega(\alpha\wedge\star_\Omega\Omega)]_{27}=[\delta_\Omega(\alpha\wedge\Omega)]_{27}
\end{array}\right.
\ee

\bigskip

\ben
d^{14}_q:\Omega_1\to\Omega_q\quad
\left\{\begin{array}{lcl}
d^{14}_1\beta & \equiv & 0\\[2pt]
d^{14}_{14}\beta & \equiv & 0\\[2pt]
d^{14}_7\beta & = & -\star_\Omega\big([d\beta]_7\wedge\Omega)\\[2pt]
d^{14}_{27}\beta & = & [d\beta]_{27} 
\end{array}\right.
\ee

\bigskip

\ben
d^{27}_q:\Omega_1\to\Omega_q\quad
\left\{\begin{array}{lcl}
d^{27}_1\gamma & \equiv & 0\\[2pt]
d^{27}_{14}\gamma & = & [\delta_\Omega\gamma]_{14}\\[2pt]
d^{27}_7\gamma & = & \star_\Omega(\delta_\Omega\gamma\wedge\star_\Omega\Omega)\\[2pt]
d^{27}_{27}\gamma & = & \star_\Omega[d\gamma]_{27}
\end{array}\right.
\ee
%

\bigskip

\noindent 
H.~{\sc Wei{\ss}}: Mathematisches Institut der Universit\"at M\"unchen, Theresienstra{\ss}e 39, D--80333 M\"unchen, F.R.G.\\
e-mail: \texttt{weiss@math.lmu.de}

\smallskip

F.~{\sc Witt}: Mathematisches Institut der Universit\"at M\"unster, Einsteinstra{\ss}e 62, D--48149 M\"unster, F.R.G.\\
e-mail: \texttt{frederik.witt@uni-muenster.de}
\end{document}